%% file: Good.tex
\documentclass[reqno]{amsart}

\usepackage{etex}
\usepackage{amsfonts}
\usepackage{amscd}
\usepackage{amsbsy}
\usepackage{mathtools}
\usepackage{amsxtra}
\usepackage{amssymb}
\usepackage{epsfig}
\usepackage{epic,eepic}
\usepackage{graphicx}
\usepackage[all]{xy}
\usepackage{xypic}
\usepackage{psfrag}
\usepackage{amsmath}
\usepackage{amsthm}
\usepackage{setspace}
\usepackage{hyperref}
\usepackage{url}
\usepackage{here}
\usepackage{todonotes}
\newlength{\halfbls}
\usepackage{xparse}
\usepackage{xspace}
\usepackage[normalem]{ulem}

\usepackage{tikz-cd}
\usepackage{tikz}
\usetikzlibrary{patterns}
\usetikzlibrary{calc}
\usetikzlibrary{decorations.pathreplacing,decorations.markings}
\usetikzlibrary{arrows}
\usetikzlibrary{decorations.markings, decorations.shapes}
\usetikzlibrary{positioning}
\usetikzlibrary{intersections}
\usetikzlibrary{decorations.pathmorphing}
\usetikzlibrary{arrows.meta,bending}

\usepackage{thmtools}

\DeclareRobustCommand{\SkipTocEntry}[9]{}

\include{macros_Good}
\usepackage[style=alphabetic,
            isbn=false,
            doi=false,
            url=false,
            backend=bibtex,
            maxnames=5
           ]{biblatex}
\defbibheading{bibliography}[\bibname]{%
  \addtocontents{toc}{\SkipTocEntry}%
  \addcontentsline{toc}{part}{References}%
  \section*{References}%
}

\bibliography{BibGood}

\usepackage{breqn}
\allowdisplaybreaks[2]

\title[The area is a good {enough} metric]
      {The area is a good {enough} metric}

\begin{document}
\author{Matteo Costantini}
\email{costanti@math.uni-bonn.de}
\address{Institut f\"ur Mathematik, Universit\"at Bonn,
Endenicher Allee 60,
53115 Bonn, Germany}
\author{Martin M\"oller}
\email{moeller@math.uni-frankfurt.de}
\author{Jonathan Zachhuber}
\email{zachhuber@math.uni-frankfurt.de}
\thanks{Research 
is partially supported  
  by the LOEWE-Schwerpunkt
``Uniformisierte Strukturen in Arithmetik und Geometrie''.}
\address{
Institut f\"ur Mathematik, Goethe--Universit\"at Frankfurt,
Robert-Mayer-Str. 6--8,
60325 Frankfurt am Main, Germany
}

\maketitle
\begin{abstract}
In the first part we extend the construction of the smooth
normal-crossing divisors  compactification  of projectivized strata of
abelian differentials given by Bainbridge, Chen, Gendron, Grushevsky
and M\"oller to the case of $k$-differentials. Since the generalized
construction is closely 
related to the original one, we mainly survey their results and
justify the details that need to be adapted in the more general context.
\par
In the second part we show that the flat area provides a
canonical hermitian metric on the tautological bundle over the projectivized
strata of finite area $k$-differentials
whose curvature form represents the first Chern class.
This result is useful in order to apply Chern-Weil theory tools.
It has already been used as an assumption in the work of Sauvaget for
abelian differentials and is also used in a paper of
Chen, M\"oller and Sauvaget for quadratic differentials.
\end{abstract}

\tableofcontents
\SaveTocDepth{1} 

\input{sec_intro} 

\input{sec_Construction}
\input{sec_good_new}





\printbibliography

\end{document}

%% file: macros_Good.tex
\DeclareMathAlphabet\mathbfcal{OMS}{cmsy}{b}{n}
\usepackage{adjustbox}

\theoremstyle{plain}
\newtheorem{Defi}{Definition}[section]  
   \newtheorem{Prop}[Defi]{Proposition}
\newtheorem{Lemma}[Defi]{Lemma}    
\newtheorem{Thm}[Defi]{Theorem}

\newtheorem{Conv}[Defi]{Convention}

\theoremstyle{remark}

\newcommand{\CC}{\mathbb{C}}

\newcommand{\NN}{\mathbb{N}}
 
\newcommand{\PP}{\mathbb{P}} 
\newcommand{\QQ}{\mathbb{Q}} 
\newcommand{\RR}{\mathbb{R}} 
 
\newcommand{\VV}{\mathbb{V}} 
\newcommand{\WW}{\mathbb{W}} 
\newcommand{\ZZ}{\mathbb{Z}}

\newcommand{\cLL}{{\mathcal L}} \newcommand{\cAA}{{\mathcal A}}   
  
 \newcommand{\cXX}{{\mathcal X}}  
 
\newcommand{\cBB}{{\mathcal B}} 
 
 \newcommand{\cOO}{{\mathcal O}}

\newcommand{\cYY}{{\mathcal Y}}

\newcommand{\bfq}{{\bf q}}

\newcommand{\bfs}{{\bf s}}
\newcommand{\bft}{{\bf t}}

\newcommand{\bfw}{{\bf w}}
\newcommand{\bfx}{{\bf x}}
\newcommand{\bfz}{{\bf z}}

\newcommand{\bfxi}{{\boldsymbol{\xi}}}

\newcommand{\bfomega}{{\boldsymbol{\omega}}}
\newcommand{\bfeta}{{\boldsymbol{\eta}}}

\newcommand{\Aut}{{\rm Aut}}

\newcommand{\area}{{\rm area}} 

\def\={\;=\;}  \def\+{\,+\,}          \def\d{\partial}

\DeclareMathAlphabet{\eucal}{U}{eus}{m}{n}
\DeclareMathAlphabet{\newcal}{U}{dutchcal}{m}{n}

      \def\wh#1{\widehat{#1}}  
\newcommand{\ol}{\overline}

\def\wP{\widehat P}  
\def\wX{\widehat X}
\def\wY{\widehat Y}
\def\wZ{\widehat Z}
\newcommand{\whmu}{\widehat{\mu}}
\newcommand{\whLa}{\widehat{\Lambda}}

       \newcommand{\ve}{{\varepsilon}}

\def\be{\begin{equation}}   \def\ee{\end{equation}}     \def\bes{\begin{equation*}}    \def\ees{\end{equation*}}
\def\ba{\be\begin{aligned}} \def\ea{\end{aligned}\ee}   \def\bas{\bes\begin{aligned}}  \def\eas{\end{aligned}\ees}

\newcommand{\omoduli}[1][g]{{\Omega \mathcal M}_{#1}}

\newcommand{\barmoduli}[1][g]{{\overline{\mathcal M}}_{#1}}

\newcommand{\Teichmuller}{Teich\-m\"uller}

\newcommand{\grc}{{\mathrm{grc}}}

\newcounter{savedtocdepth}
\newcommand*{\SaveTocDepth}[1]{%
  \addtocontents{toc}{%
    \protect\setcounter{savedtocdepth}{\protect\value{tocdepth}}%
    \protect\setcounter{tocdepth}{#1}%
  }%
}

\DeclareDocumentCommand{\LMS}{ O{\mu} O{g,n}} {\Xi\overline{\mathcal{M}}_{#2}(#1)}
\DeclareDocumentCommand{\kLMS}{ O{\mu} O{g,n}} {\Xi^k\!\overline{\mathcal{M}}_{#2}(#1)}

\newcommand*{\mskd}{multi-scale $k$-differential\xspace}

\newcommand*{\mskds}{multi-scale $k$-differentials\xspace}

\newcommand{\komoduli}[1][g,n]{{\Omega^k\!\mathcal M}_{#1}}

\newcommand{\PPer}{{\rm PPer}}
\newcommand{\tkd}[1][\wh{\Gamma}]{\mathfrak{W}^{k}(#1)}

\newcommand{\tkdnoab}[1][\wh{\Gamma}]{\mathfrak{W}^{k}_{\operatorname{na}}(#1)}
\newcommand{\tkdpm}[1][\wh{\Gamma}]{\mathfrak{W}^{k}_{\operatorname{pm}}(#1)}
\newcommand{\tkdpms}[1][\wh{\Gamma}]{\mathfrak{W}^{k,s}_{\operatorname{pm}}(#1)}

\newcommand*{\Tw}[1][\wh\Gamma]{\mathrm{Tw}_{#1}}  
\newcommand*{\sTw}[1][\wh\Gamma]{\mathrm{Tw}_{#1}^s}  
\newcommand*{\LRT}[1][\wh\Gamma]{\mathrm{T}_{#1}}  
\newcommand*{\sLRT}[1][\wh\Gamma]{\mathrm{T}_{#1}^s}  

\newcommand{\prodt}[1][\lceil i \rceil]{t_{#1}}

\newcommand{\prodttop}[1][\lceil e^+\rceil]{t_{#1}}
\newcommand{\prodtbot}[1][\lceil e^-\rceil]{t_{#1}}

\newcommand{\hslashslash}{%
\lapbox[\width]{-0.15em}{\raisebox{.05ex}{%
    \scalebox{.7}{%
      \rotatebox[origin=c]{22}{$-$}%
    }%
  }%
}}

\newcommand{\bslash}{%
  {%
   \vphantom{d}%
   \ooalign{\kern.05em\smash{\hslashslash}\hidewidth\cr$b$\cr}%
   \kern.05em
  }%
}

\newcommand{\tui}[1][i]{t_{\lceil #1 \rceil}}

%% file: sec_intro.tex
\section{Introduction}
\label{sec:intro}

A flat surface $(X,\omega)$ is a Riemann surface together with
a non-zero holomorphic one-form. Interest in flat surfaces stems 
from dynamics of polygonal billiards and this paper contributes
to justifying the foundations for the efficient computation of
invariants like Siegel-Veech constants of these billiards.
A natural invariant of a flat surface is the flat
area ${\rm vol}(X,\omega)$, the area taken with respect to the
form~$|\omega|$. As such, it defines a hermitian metric~$h$ on the tautological
line bundle $\cOO(-1)$ over the projectivized strata $\PP\omoduli(\mu)$,
the moduli space parameterizing flat surfaces whose zeros and poles are
of a fixed type~$\mu = (m_1,\ldots,m_n)$. This metric does not extend smoothly
over the boundary, as the area of a flat surface tends to~$\infty$
when~$X$ acquires an infinite flat cylinder, i.e.\ when~$\omega$ acquires
a simple pole.  In Chern-Weil theory applications, it suffices
to show that the curvature form of the metric connection associated
to the metric~$h$  represents the first Chern
class of~$\cOO(-1)$ on a suitable compactification. This
has been used as assumption by Sauvaget in~\cite{SauvagetMinimal} for
Masur-Veech volumes of the minimal strata of abelian differentials. While
a workaround for this has been given in~\cite{CMSZ}, the computation
of the volume of individual spin components in loc.\ cit.\ is still
based on that assumption. Moreover, {the} 
paper~\cite{CMS2} extends this line of thought to quadratic differentials.
There, too, the volume of the canonical
double cover (see Section~\ref{sec:period}) provides a natural hermitian metric.
Even for principal strata, where the Hodge bundle provides a smooth
compactification, we do not see an easy route to prove the claim
in the title, see the subtleties explained below.  This paper consequently
makes full use of the smooth compactification of strata of abelian differentials
constructed in \cite{LMS}. {Yet another application is a growth justification
in the recent computation of the volume of moduli spaces of flat surfaces
(in the sense of Veech (\cite{VeechFlat}) by Sauvaget \cite{SauvagetFlat}).}
\par
\medskip
Given the applications in mind, the first part of this paper is a survey
about the construction of the smooth compactification and  the formal
justification of the tempting claim that the construction extends
to $k$-differentials, if the notions are appropriately
adapted in the same way as \cite{kdiff} adapts \cite{strata}. 
\par
\medskip
\paragraph{\bf The compactification}
Let $\mu = (m_1,\ldots,m_n)\in \ZZ^n$ be a type of a meromorphic $k$-differential,
i.e. $m_i$ are integers such that  $\sum m_i = k(2g-2)$. Let $\komoduli(\mu)$ 
be the moduli space of triples $(X,\bfz,q)$ consisting of a smooth curve~$X$
of genus~$g$ with marked points $\bfz = (z_1,\ldots,z_n)$ and a $k$-differential
having zeros or poles of order~$m_i$ at the points~$z_i$. We summarize the properties
of our compactification $\PP \kLMS$ of this moduli space of $k$-differentials.
The canonical cover construction and related notions are recalled in
Section~\ref{sec:period}.
\par
\begin{Thm} \label{thm:kLMS}
There exists a complex orbifold $\kLMS$, the moduli space of \mskds,
with the following properties.
\begin{itemize}
\item[i)] The space $\komoduli(\mu)$ is dense in $\kLMS$.
\item[ii)] The boundary $D = \kLMS \smallsetminus \komoduli(\mu)$ is a
normal crossing divisor.
\item[iii)] The rescaling action of $\CC^*$ on $\komoduli$ extends to $\kLMS$
and the resulting projectivization $\PP \kLMS$ is a compactification
of $\PP \komoduli(\mu)$.
\item[iv)] Via the canonical cover construction, the space $\kLMS$ has
a map to the compactification $\LMS[\wh{\mu}][\wh{g},\{\wh{n}\}]$ of
the corresponding stratum of abelian differentials with partially labeled
points. In the interior, this map is a closed immersion.
\end{itemize}
\end{Thm}
\par
Here we only prove that $\kLMS$ is a 'moduli space' in a very weak form,
namely by exhibiting what its complex points correspond to, the
\mskds introduced below. We leave it to the interested reader to adapt
the functor from \cite{LMS} to the context of $k$-differentials. Recall
that a $k$-differential is called \emph{primitive} if it cannot be written
as $d$-th power of a $k/d$-differential for any $d>1$. In general our
notion of $k$-differentials does not imply primitivity. This is convenient
for defining twisted $k$-differentials below, but as a consequence
the spaces in Theorem~\ref{thm:kLMS} have many connected components,
those consisting of $k$-th powers of abelian differentials having
dimension one more than the other components.
\par
Besides the normal crossing boundary, the most relevant property for us is the
existence of a convenient coordinate system, given by {\em perturbed
period coordinates}. To introduce this, we first have to explain
how to parameterize boundary points of $\kLMS$.
\par
Let $\Gamma = (V,H,E,g)$ be a stable graph, as in \cite{acgh2}, where $V$ is
the set of vertices, $H$ is the set of half-edges, $E$ is the set of edges
and  $g$ is the genus assignment.
A {\em level graph} is a stable graph together with a weak total order on
the set of vertices, which is determined by a level function normalized to take
values in $0,-1,\ldots,-L$, with zero being the top level and the rest often
referred to as \emph{lower levels}. An edge of a level
graph is called \emph{horizontal} if  it is adjacent vertices are
on the same level
and \emph{vertical otherwise}. This leads to a partition of the edges
$E = E^h \cup E^v$ and we use the adjectives horizontal and vertical for the
corresponding nodes accordingly. 
An {\em enhanced level graph} is a level graph together with an enhancement
$\kappa \colon H \to \ZZ$ on the half-edges that specifies the number of
prongs of the differential at the corresponding marked point, see
Section~\ref{sec:period} for the full definition.
\par
Each of the levels of~$\Gamma$ thus specifies a moduli space of
$k$-differentials, the type being given by the enhancement and the restrction
of~$\mu$ to the legs at that level. A collection of these differentials, one
for each level, on a given pointed stable curve $(X,\bfz)$ is called
{\em twisted differential} and we call a twisted differential {\em compatible
with the enhanced level graph~$\Gamma$} if the underlying graph of~$\Gamma$
is the stable graph of $(X,\bfz)$ and if the collection moreover satisfies
the global $k$-residue-condition (GRC) from \cite{kdiff}.
A {\em \mskd} is a twisted differential compatible with~$\Gamma$
up to projectivization of the  levels below zero, together with the choice
of an equivalence class of prong-matchings. The details are given
in Section~\ref{sec:Constr} using the notion of level rotation torus.
Leaving them aside, we can now describe the coordinates. 
\par
\begin{Prop} \label{prop:PPerIntro}
In a neighborhood~$U \subset \kLMS$ of every point in the boundary stratum
corresponding to an enhanced level graph~$\Gamma$ with~$L+1$ levels and
$h$ horizontal edges there
is an orbifold chart given by the perturbed period map
\bes \PPer : U \to \CC^{|E^h|} \times \CC^{L+1} \times \prod_{i=0}^{L}
\CC^{\dim E_{(-i)}^{\grc} -1 }\,,
\ees
where $E_{(-i)}^{\grc}$ is some eigenspace in homology constrained by the
GRC and where the corresponding coordinates are obtained by integrating
perturbations of the twisted differential against these homology classes.
\end{Prop}
\par
In this proposition, the first set of coordinates in~$\CC^{|E^h|}$ measures the
opening of horizontal nodes and the second set in~$\CC^{L}$ measures the
rescaling of the differentials on each level. Neither of them
is a period, in fact they are exponentials, respectively roots, of periods. 
The statement about integration is intentionally vague, since we are not
exactly integrating the (roots of) $k$-differentials parameterized by~$U$,
but its sum with a modification differential, as constructed in
Section~\ref{sec:Constr}. Moreover, the path of integration is not between
the zeros of those differentials but between neighboring points, thus the
name 'perturbed'. Technically important is that these perturbations go
to zero faster than the rescaling of the $k$-differential. The map
$\PPer$ depends on many choices (see Section \ref{sec:PPer} for more details),
however they are irrelevant for
many local computations.
\par
\medskip
\paragraph{\bf Boundary divisors}
To a first approximation the boundary divisors, i.e., the irreducible components
of the boundary  $\kLMS \smallsetminus \komoduli(\mu)$, are
given by graphs with one level and a single horizontal edge, and by graphs
with two levels and no horizontal edge. However, in the setting
of $k$-differentials the level graph does not
specify the boundary divisor uniquely. In Section~\ref{sec:period} we recall
the notion of canonical $k$-cover, which is unique for $k$-differentials
on smooth curves, but not in the stable case. An example for two
different covers that give rise to different components of the boundary
is given by \cite[Figure~2]{kdiff}. In fact, the residue conditions are
different in the two cases. Consequently, as second approximation
the choice of a cyclic $k$-cover $\pi \colon \wh\Gamma \to \Gamma$ compatible
with the canonical covers of the components (see Section~\ref{sec:period}
for the definition of both notions) characterizes boundary components.
\par
\begin{Prop}
For each $k$-cyclic cover $\pi \colon \wh\Gamma \to \Gamma$ of
enhanced level graphs with~$\Gamma$ of type~$(g,n)$ there is
a boundary stratum $D_{\wh\Gamma}$ of  the compactification $\kLMS$.
Each $D_{\wh\Gamma}$ is commensurable to the product of the  level-wise
projectivised moduli space of twisted differentials on $\wh\Gamma$.
\end{Prop}
\par
Here 'commensurable' is a shorthand for the existence of a complex space
with a finite map to the two spaces in question.
\par
We will not address the subtle question of connectivity of those~$D_{\wh\Gamma}$.
The details of the construction of a space that admits a finite covering
to both~$D_{\wh\Gamma}$ and the product level-wise projectivised moduli spaces
is given in \cite[Section~4.2]{StrataEC}. There, the construction is given
for Abelian differentials, but it can verbatim be applied for $k$-differentials,
too.
\par
\medskip
\paragraph{\bf The metric} We now return to our primary goal. The statement
is about flat surfaces of finite area, so we suppose from now on
that $m_i >-k$. If
$\pi \colon \wX \to X$ denotes the canonical covering associated
with $(X,q) \in \komoduli(\mu)$ such that $\pi^*q = \omega^k$ is a
$k$-th power, then the definition 
\be \label{eq:defh}
h(X,q)^{1/k} \= \area_{\wX}(\omega) \= \frac{i}{2} \int_{\wX} \omega \wedge 
\overline{\omega}
\ee
provides the tautological bundle $\cOO(-1)$
on $\PP \komoduli(\mu)$ with a hermitian metric~$h$.
The moduli space $\PP \komoduli(\mu)$ has, besides the nice 
compactification $\kLMS$ discussed above, a highly singular
compactification, the {\em incidence variety compactification}
$\PP \overline{\Omega^k\mathcal{M}}_{g,n}(\mu)$ that has been studied in \cite{strata}
and \cite{kdiff}. It is the closure of $\PP \komoduli(\mu)$ inside the
projectivized bundle of $k$-fold stable differentials twisted by
the polar part of~$\mu$. This projectivized bundle provides
an extension of the tautological bundle $\cOO(-1)$, whose restriction
to the incidence variety compactification we denote by the same symbol.
\par
There is a natural forgetful map  $\varphi \colon \PP \kLMS \to \PP
 \overline{\Omega^k\mathcal{M}}_{g,n}(\mu)$, which is an isomorphism restricted to
$\PP \komoduli(\mu)$. The pullback of $\cOO(-1)$ thus provides an extension
of the tautological bundle on $ \PP \kLMS$ that we still denote by
the same symbol. It is this bundle whose Chern classes are relevant
(\cite{SauvagetMinimal}, \cite{CMSZ}) for computation of Masur-Veech volumes
and Siegel-Veech constants. Our main theorem is:
\par
\begin{Thm} \label{thm:areagood}
{The curvature form  $\tfrac{i}{2\pi}[F_h]$ of the metric~$h$
is a closed current on $\PP \kLMS$ that represents the first Chern class
$c_1(\cOO(-1))$. More generally, the $d$-th wedge power of the curvature
form represents $c_1(\cOO(-1))^d$ for any $d \geq 1$.}
\end{Thm}
\par
{In an earlier version of the paper we had claimed that the metric~$h$ is
good in the sense of Mumford. This is not true at boundary points where
there are both horizontal and vertical edges, as explained in
Section~\ref{sec:good}. We thank Duc-Manh Nguyen for bringing this to
our attention.}
\par
{ More precisely, in} the case of only horizontal nodes the metric
diverges as we approach the boundary. However in perturbed period coordinates
coordinates the local calculation is essentially the calculation of Mumford
for the special case of elliptic curves (times the number of horizontal
nodes).
\par
In the absence of horizontal nodes, the metric  extends continuously over
the boundary. This fits with the intuition that the area of the lower
level surfaces goes to zero. The area is not a $C^2$-function near vertical
notes, but the second derivative is integrable, which turns out to be
good enough. In the presence of both horizontal
and vertical edges we estimate directly the growth of the curvature form
to justify Theorem~\ref{thm:areagood}. This is a delicate computation
that makes full use of the coordinate system that we have near the
boundary, which in turn is the main reason for working with the
compactification we constructed rather than say simple with the
full Hodge bundle $\PP  \overline{\Omega^k\mathcal{M}}_{g,n}$.
\par

\subsection*{Acknowledgements}

We are very grateful to the Mathematisches Forschungsinstitut Oberwolfach
for providing a stimulating atmosphere and would like to thank Matt Bainbridge,
Dawei Chen, Quentin Gendron, Sam Grushevsky, {Duc-Manh Nguyen} and Adrien Sauvaget
for inspiring discussions. We thank an anonymous referee for valuable
comments and corrections.

%% file: sec_Construction.tex
\section{Period coordinates and canonical covers of $k$-differentials}
\label{sec:period}

In this section we summarize well-known results about period
coordinates, but also recall the period coordinates along the
boundary strata of the incidence variety compactification from \cite{kdiff}.
We start by recalling properties of the canonical $k$-cover.
\par
Let $X$ be a Riemann surface and let~$q$ be a meromorphic $k$-differential
of type~$\mu$. This datum defines (see e.g.\ \cite[Section~2.1]{kdiff})
a $k$-fold cover $\pi \colon \wX \to X$ such that $\pi^* q = \omega^k$
is the $k$-power of an abelian differential. Note that~$\wX$ is
disconnected, if $q$ is a $d$-th power of a $k/d$-differential for some~$d>1$.
This differential~$\omega$ is of type 
$$
\whmu \,:=\, \Bigl(\underbrace{\wh m_1, \ldots, \wh m_1}_{\gcd(k,m_{1})},\,
\underbrace{\wh m_2,  \ldots, \wh m_2}_{\gcd(k,m_{2})} ,\ldots,\,
  \underbrace{\wh m_n, \ldots, \wh m_n}_{\gcd(k,m_{n})} \Bigr)\,,
$$
where $\wh m_i := \tfrac{k+m_{i}}{\gcd(k,m_{i})}-1$. We let
$\wh{g} = g(\wX)$ and $\wh{n} = \sum_i \gcd(k,m_{i})$. The type
of the covering determines a natural subgroup $S_{\whmu} \subset S_{\wh{n}}$
of the symmetric group that allows only the permutations of each the
$\gcd(k,m_i)$ points corresponding to a preimage of the $i$-th point.
\par
We fix {\em once and for all a primitive $k$-th root of unity~$\zeta$}. The
Deck group of~$\pi$ contains a unique element~$\tau$ such that
$\tau^* \omega = \zeta \omega$. We fix this automorphism as well.
We denote by $\bfz = (z_1,\ldots,z_n)$ the tuple of marked points in~$X$.
The preimages in~$\wX$ of these marked points give a tuple that is
labeled up to the action of~$S_{\whmu}$ and which we denote by~$\wh{\bfz}$.
By the canonical cover construction there is an isomorphism of orbifolds
between the moduli space $\komoduli(\mu)$ and the space of
\be \label{eq:point_Hk}
\{(\wh{X},\wh{\bfz},\omega, \langle \tau \rangle) \,:\, \tau \in \Aut(\wh{X}),
\quad \mathrm{ord}(\tau)\=k, \quad
\tau^* \omega \= \zeta_k \omega \}
\ee
which has a natural closed immersion to $\omoduli[\wh{g},\wh{n}](\wh{\mu})
/S_{\whmu}$.
\par
For the analogous statements about coverings in the stable case we first
need to define twisted $k$-differentials and further preparation.
An {\em enhanced level graph} for $k$-differentials is a level graph together
with an enhancement map $\kappa: H \to \ZZ$
on the half-edges, satisfying  the following properties:
\begin{itemize}
\item[i)] If $h$ and $h'$ are paired to an edge, then $\kappa(h) +
\kappa(h') = 0$.
\item[ii)] At a leg $h \in H \smallsetminus E$ with order $m_i$,
we impose that $\kappa(h) = m_i + k$.
\item[iii)] At each vertex $v \in V(\Gamma)$ 
\bes
k(2g(v) -2)  \= \sum_{h \vdash v} (\kappa(h) - k)\,.
\ees
\end{itemize}
\par
The next notion is a combinatorial model of canonical curves, as they
occur in the limit when tending to stable curves. Let $\wh\Gamma$ be
an enhanced level graph for abelian differentials of type $(\wh{g},\wh{n})$
and let $\Gamma$ be an enhanced level graph for $k$-differentials
of type~$(g,n)$. A \emph{(cyclic) $k$-cover of nhanced level graphs}
$\pi \colon \wh\Gamma \to \Gamma$ is a morphism f graphs (with legs),
given as the quotient map by a graph automorphism $\tau \in \Aut(\wh\Gamma)$
of order~$k$ that preserves levels, the orders~$\wh{m}_i$ and enhancements,
and with the following two properties: An edge~$e$ has $\gcd(\kappa_e,k)$
preimages and a marked point of type~$m_i$ has $\gcd(m_i,k)$ preimages.
\par
We next give a the definition of a twisted differential. The case
of $k$-differentials is reduced to the case of abelian differentials.
Let $\wh{\Gamma}$ be an enhanced level graph for the stable
curve $(\wh{X},\wh\bfz)$. A \emph{twisted $1$-differential~$\bfomega =
  (\omega_v)$ compatible with $\wh\Gamma$} is a collection
of differentials on $\wh{X}_v$ for each vertex $v \in V(\wh{\Gamma})$ of
type given by the markings and enhancements, i.e., $\mathrm{ord}_e(\omega_v)
= \kappa_e -1$ for
each edge~$e$ adjacent to~$v$. This collection is required to satisfy the usual
residue condition at horizontal nodes and moreover the global residue
condition (GRC), see~\cite{kdiff} for details on this. Later it will be
convenient to group together the differentials on all vertices of the
same level and we thus write $\bfomega = (\omega_i)_{i \in L(\wh\Gamma)}$.
\par
A collection of $k$-differentials on~$X$ naturally defines a $k$-cover for
each component of~$X$, but it does not uniquely define how to glue them to
a stable curve~$\wh{X}$. A $k$-cover of level graphs contains this gluing
information. Let ${\Gamma}$ be an enhanced level graph for the stable
curve $({X},\bfz)$. A \emph{twisted $k$-differential compatible with $\Gamma$} is a
collection $\bfq = (q_v)$ of differentials on $X_v$ for each vertex $v \in V({\Gamma})$
such that the pullback to the stable curve $\wh{X}$ given by the canonical
covers induced by~$q_v$ and some $k$-cover $\pi \colon \wh\Gamma \to \Gamma$
is a $k$-th power of a twisted abelian differential. (A formulation of
this condition directly on~$X$ is given in \cite{kdiff}. This includes
however a quite complicated formulation of the GRC.) As above we group
these differentials according to levels and write $\bfq = (q_i)_{i \in L(\Gamma)}$.
We say that $\bfq$ is \emph{compatible with $\wh\Gamma$}, if the above condition
holds for a chosen $k$-cover $\pi \colon \wh\Gamma \to \Gamma$. By definition
we can specify a twisted $k$-differential either by $(X,\bfz,\bfq, \pi)$
or by $(\wh{X},\wh{\bfz}, \bfomega,\tau)$, where $\tau$ is an automorphism
of order~$k$. 
\par
The starting point for the construction of the compactification is the
moduli space of twisted $k$-differentials compatible with~$\wh\Gamma$,
which we denote by $\tkd$ suppressing the dependence on the initial
type~$\mu$. By definition this is a subspace defined by the GRC inside
a product of strata. The main result of \cite{kdiff} implies that points
in $\tkd$ are smoothable: they arise as limits of $k$-differentials
in $\komoduli(\mu)$. 
\par
\bigskip
Next we define the subspaces of homology that we use for period
coordinates. We fix some reference smooth surface $\Sigma$ of genus~$g$
with $n$ marked points that we partition as $P \cup Z$ according
to the poles of order $\leq -k$ among~$\mu$ and the 'zeros', i.e.
points with order~$>-k$. We let $\wh{\Sigma}$ be a model for the canonical
covering surface, which is of genus $\wh{g}$, and which comes with a
map $\pi \colon \wh{\Sigma} \to \Sigma$. We let $\wh{P}$ and $\wh{Z}$
be the preimages of~$P$ and $Z$. They now correspond indeed
to the zeros and poles of the type~$\wh{\mu}$.
\par
\begin{center}
	\begin{tikzcd}
\wh{P}, \wh{Z} \arrow{d}{\pi}  & \wh{\Sigma}  \arrow{d}{\pi} &
\wh{\Lambda}^\circ \supset \wh{\Lambda}
\arrow{d}{\pi} 
\\
P, Z  & \Sigma  & \Lambda^\circ \supset \Lambda   
	\end{tikzcd}
\end{center}
If $X$ is a stable curve and $\pi$ a covering as above,
we may find a multicurve $\wh{\Lambda}$ in~$\wh{\Sigma}$ mapping under~$\pi$
to the multicurve~$\Lambda$ in~$\Sigma$, such that~$\wX$ and~$X$
are obtained by pinching $\wh{\Sigma}$ and $\Sigma$ along  $\wh{\Lambda}$
and $\Lambda$ respectively.
\par
Recall (\cite[Section~2]{kdiff}) that the moduli space $\komoduli(\mu)$
of $k$-differentials is locally modeled on the $\omega$-periods of the eigenspace
\bes
E(\wh{\Sigma} \smallsetminus \wh{P}, \wh{Z}) \=
H_1(\wh{\Sigma} \smallsetminus \wh{P}, \wh{Z},\CC)_{\tau = \zeta}\,.
\ees
Similarly, we can describe local coordinates for the components
of a twisted $k$-differential on a stable curve~$X$ with enhanced level
graph~$\Gamma$ (not yet imposing full compatibility, i.e.\ the GRC).
Let $\Lambda^\circ$ be an open thickening of $\Lambda$. We
let $\Lambda^\pm$ be the upper and lower boundaries of $\Lambda^\circ$.
The level structure on $\Gamma$ organizes $\Sigma \smallsetminus \Lambda$
into levels $\Sigma_{(i)}$ and we denote the adjacent poles, zeros
and boundaries $\Lambda^\pm$ with the subscript~$(i)$. All the
notation is applied with a hat to the corresponding objects on the $k$-cover.
The level-$i$ component of the twisted differential is thus
modeled on
\be \label{eq:H1leveli}
E_{(i)} \= H_1(\wh\Sigma_{(i)} \setminus \{\wP_{(i)} \cup \whLa_{(i)}^\circ \},
\whLa_{(i)}^+ \cup \wZ_{(i)},\CC)_{\tau = \zeta}\,.
\ee
We can now restate the main dimension estimate in the proof of
\cite[Theorem~6.2]{kdiff}.
\par
\begin{Prop} \label{prop:dimper} The moduli space of twisted $k$-differentials
compatible with an enhanced level graph~$\Gamma$ is locally modeled on
the $\omega_{(-i)}$-periods of $\prod_{i=0}^{L} E_{(-i)}^{\grc}$, where $E_{(i)}^{\grc}
\subseteq E_{(i)}$ is the subspace at level~$i \in L(\Gamma)$ cut out by the
global residue condition. The dimensions of these subspaces is given by
\bes
\sum_{i=0}^{L} \dim_{\CC}E_{(-i)}^{\grc} \= \dim_{\CC} \komoduli\,-\,|E^h|\,,
\ees
where $E^h$ is the set of horizontal edges of~$\Gamma$.
\end{Prop}

\section{Construction of $\kLMS$}
\label{sec:Constr}

In this section we recall the main technical tools from
\cite{LMS}, construct the compactification and eventually
prove Theorem~\ref{thm:kLMS}. The definitions in
Section~\ref{sec:DegUndeg}--Section~\ref{sec:TopoComp} are direct
adaptations of the abelian case by working on the canonical $k$-covers.
Avoiding the discussion of \Teichmuller\ spaces means omission of
that aspect but also a simplification of notations. In the remaining
sections we have to ensure at some places that constructions can
be performed $\tau$-equivariantly.

\subsection{Degeneration, undegeneration}
\label{sec:DegUndeg}

We describe here two types of maps between level graphs~$\Gamma$
that encode the degeneration of curves, together with the compatible maps
between the coverings graphs~$\wh\Gamma$ that form part of the
degeneration datum. In fact, it is easier
to first describe the inverse process of undegeneration that
encodes all the $k$-differentials in a neighborhood of a given
degenerate $k$-differential.
\par
Let $\pi \colon \wh\Gamma \to \Gamma$ by a cyclic $k$-covering
of enhanced level graphs with~$L+1$ normalized levels. For any subset
$I \subset \{1,\ldots,L\}$,
to be memorized as the {\em the
level passages that remain}, we define the {\em vertical undegeneration}
$\delta_I$ as the following contraction of certain vertical edges.
An edge~$e$ is contracted by $\delta_I$ if and only it crosses the level
passages indexed by the complement of $I$. Vertices are merged if an edge
connecting them has been contracted. This edge
contraction is performed simultaneously on the domain and range of~$\pi$
and induces a cyclic $k$-covering $\delta_I(\wh\Gamma) \to \delta_I(\Gamma)$
that we abbreviate as $\delta_I(\pi)$. We write $\delta_I(j)$ for the
image of the $j$-th level under~$\delta_I$.
\par
Moreover, we define for any subset $E_0 \subset E^{\rm h}$ of the
horizontal edges of~$\Gamma$ the {\em horizontal undegeneration} $\delta_{E_0}$
to be the edge contraction that contracts precisely the edges in $E_0$
in $\Gamma$. Contracting simultaneously also on the $\pi$-preimages of~$E_0$
in $\wh{\Gamma}$, we obtain a new cyclic $k$-covering $\delta_{E_0}(\pi):
\delta_{E_0}(\wh{\Gamma})\to \delta_{E_0}(\Gamma)$.
\par
A general undegeneration is a pair $\delta = (\delta_I, \delta_{E_0})$,
defined by composing  a horizontal and a vertical  undegeneration in either
order. A {\em degeneration} is the inverse of an undegeneration.
We write $\wh\Gamma' \rightsquigarrow \wh\Gamma$ for a general
degeneration of level graphs and~$\delta^{\rm ver}$ and~$\delta^{\rm hor}$
for the two constituents of an undegeneration~$\delta$.
\par
For any edge $e$ of $\wh\Gamma$, we will denote by  $e^-$  both the point
in $\wX$ at the bottom end of the edge~$e$ and its level. The meaning should
be clear from the context. Similarly, $e^+$ refers to the top end.

\subsection{Prong-matchings as extra structure on twisted differentials}
\label{sec:PM}

We start with the definition of prong-matchings and the welded surfaces
constructed from these. Given a differential~$\omega$ on~$X$ that
has been put in standard form, which is $z^\kappa dz/z$ if $\kappa > 0$ or
$(z^\kappa + r)dz/z$ if $\kappa <0$, the {\em prongs} are the~$|\kappa|$
tangent vectors $e^{2\pi i j/|\kappa|} \tfrac{\partial}{\partial z}$ for
$\kappa >0$ and $-e^{2\pi i j/|\kappa|} \tfrac{\partial}{\partial z}$
for $j=0,\dots,|\kappa|-1$. At simple poles (i.e.\ for $\kappa=0$),
prongs are not defined.
\par
We now get back to a twisted $1$-differential~$(\wX,\bfomega,\wh\Gamma)$. Define
a {\em local prong-matching $\sigma_e$} at the vertical edge~$e$ of~$\wh\Gamma$
to be a cyclic order-reversing bijection between the $\kappa_e$ prongs of
$\bfomega$ at the upper and lower end of~$e$. A {\em global prong-matching} is
a collection $\sigma = (\sigma_e)_{e \in E(\wh\Gamma)}$ of local prong-matchings.
If the twisted $1$-differential stems from a twisted $k$-differential
$(X,\bfz,\bfq, \pi)$, i.e. if it contains the additional information of
the automorphism~$\tau$, we require moreover that the global prong-matching
is equivariant with respect to the action of~$\tau$ permuting the edges
and multiplying the local coordinates~$z$ by~$\zeta$. 
\par
A global prong-matching~$\sigma$ on~$\wX$ gives an {\em almost-smooth
surface}~${\wX}_\sigma$, i.e.\ a smooth surface except for nodes
corresponding to the horizontal nodes of~$\wX$, constructed by the
following procedure of {\em welding}. Take the partial normalization
of~$\wX$ separating branches at vertical nodes and
perform the real oriented blowup of each pair of preimages. Then identify
the boundary circles isometrically so as to identify
boundary points that are paired by the prong-matching. More details of
the construction can be found in Section~4 of \cite{LMS}, see also
\cite{acgh2}.  (We only use subscripts~$\sigma$ to denote weldings here and
suppress the overline used in loc.~cit.\ to avoid double decorations.)
Horizontal nodes remain untouched in the welding procedure.
\par
On almost-smooth surfaces any {\em good arc~$\gamma$}, i.e.\ any arc transversal
to the seams created by welding, has a well-defined {\em turning number}
with respect to the flat structure~$\omega$, that we denote by~$\rho(\gamma)$.
\par
\medskip
Adding the information of prong-matching to points in the space of twisted
differentials $\tkd$ will give us a finite covering space as follows.
\par
We define the set $\tkdpm$ to be tuples $(\wh{X},\wh{\bfz}, \bfomega, \tau, 
\wh{\sigma})$ consisting of a point  $(\wh{X},\wh{\bfz}, \bfomega, \tau)
\in \tkd$ together with a prong-matching~$\wh{\sigma}$. There is an obvious
notion of parallel transport of prong-matchings (since the finitely many tangent vectors
$\pm e^{2\pi i j/|\kappa|} \tfrac{\partial}{\partial z}$ depend continuously
on the twisted differential) that allows to lift inclusions of contractible
open sets $U \to \tkdnoab$ uniquely to maps $U \to \tkdpm$. Requiring that
these lifts are holomorphic local homeomorphism provides $\tkdpm$ with a
complex structure so that $\tkdpm \to \tkdnoab$ is a covering map.

\subsection{The level rotation torus}
\label{sec:LRT}

\par
Our compactification combines the geometry of moduli spaces of
$k$-differentials of lower complexity and aspects of a toroidal
compactification. The torus action for the latter is given by the level
rotation torus that we now define. To describe various group actions on
prong-matchings, we view $\wh\Gamma$ as a graph with~$L$ level~{\em passages},
the first from level~$0$ to level~$-1$, the second from level~$-1$ to
level~$-2$ etc. This is summarized by:
\par
\begin{Conv}
Levels are indexed by negative integers $0,-1,\dots,-L$, while level
passages are indexed by positive integers $1,\dots,L$.
\end{Conv}
\par
The unit vector $e_i$ in the {\em level rotation group} $R_{\wh\Gamma} \cong
\ZZ^{L}$ acts on the set of prong-matchings by shifting the prong-matching
for each edge crossing the $i$-th level passage by one counterclockwise turn.
Of particular importance
is the subgroup $\Tw$ of $R_{\wh\Gamma}$ that fixes all prongs, the {\em twist
group}. The {\em (reduced) level rotation torus} $\LRT$ is the quotient
\bes
\LRT \= \CC^{L}/\Tw\,.
\ees
(Here reduced refers to the fact that $\LRT$ does not rotate the
top level. We will introduce this action separately for projectivization
and usually drop the adjective 'reduced'.)  The level rotation torus
can also be characterized (\cite[Proposition~5.4]{LMS}) by its joint
action on edge and level parameters, namely as the connected
component of the identity of the subtorus
\be
\Bigl\{ \left( (r_i, \rho_e)\right)_{i,e} \in
(\CC^*)^{L} \times (\CC^*)^{E(\wh\Gamma)} \,| \,
  r_{|e^-|} \dots  r_{|e^+|+1} = \rho_e^{\kappa_e} \,\,\text{for all~$e \in E(\wh\Gamma)$}\Bigr\}.
\ee
This characterization makes the reason for introducing~$\LRT$ apparent,
as there is a natural action of the level rotation torus on
$\tkdpm$ given by
\ba\label{eq:actionLRT}
\LRT \times \tkdpm \quad &\to \quad \tkdpm \\
(r_{|i|}, \rho_e)\,\ast\, \bigl(\wX,(\omega_{(i)}),\,  (\sigma_e)\bigr)
&\=  \biggl(\wX, \bigl(r_{|i|}\dots r_{1}\omega_{(i)}\bigr),\,
(\rho_e \ast \sigma_e)\biggr)
\ea
where $\rho_e \ast \sigma_e$ is the prong-matching~$\sigma_e$
post-composed with the rotation by~$\arg(\rho_e)$ (if the full Dehn twist
around~$e$ corresponds to angle~$2\pi$, equivalently by the rotation
by $\kappa \arg(\rho_e)$ for the angle in the flat metric). We alert the
reader that this action uses the 'triangular' basis, where the $i$-th
component of $\LRT$ rotates the $i$-th level and all level below it by
the amount~$r_{|i|}$. 
\par
To obtain orbifold charts we need to define  the {\em simple
twist group} $\sTw \subseteq \Tw$ as the twist group elements generated
by rotations of one level at a time, i.e.\
\bes
\sTw \= \oplus_{i=1}^{L} \Tw[\delta_{i}(\wh\Gamma)]\,.
\ees
We can now define the {\em (reduced) simple level rotation torus} as
\be
\sLRT \= \CC^{L} /\sTw\,.
\ee
In order to describe the action of these tori we will need
the integers
\be \label{eq:defai}
\ell_i \= {\rm lcm}_{e \in E(\delta_{i}(\wh\Gamma))} k_e\,,\quad
\text{and} \quad m_{e,i} \= \ell_i/\kappa_e
\ee
for $i=1,\ldots,L$ and $e \in E(\wh\Gamma)$. Now Proposition~5.4
in loc.\ cit.\ moreover shows that there is an identification
$\sLRT \cong (\CC^*)^{L}$ such that the projection $\sLRT \to \LRT$
is given in coordinates by
\be \label{eq:sLRTpara}
(t_i) \,\mapsto\,
(r_i, \rho_e) \= \Bigl(t_i^{\ell_i}\,,\, \prod_{i=|e^-|}^{|e^+|+1}
t_{i}^{\ell_i/\kappa_e} \Bigr)
\ee
\par
The composition of this parametrization~\eqref{eq:sLRTpara} of
$\LRT$ by $\sLRT$ with the action~\eqref{eq:actionLRT} gives an action
of $\bft = (t_i) \in \sLRT$ on welded surfaces and we denote the
image of~$\wX_\sigma$ under the action of $\bfs$ by $\wX_{\bfs \cdot \sigma}$.
\par

\subsection{The compactification as topological space}
\label{sec:TopoComp}

We start with the definition of $\kLMS$ as a set. For each $k$-cyclic covering
$\pi:\wh{\Gamma}\to \Gamma$ 
we define the {\em boundary stratum} $\Omega^k\cBB_{\wh\Gamma} = \tkdpm / \LRT$
and we define the set
\be
\kLMS \= \coprod_{\pi:\wh\Gamma \to \Gamma}\Omega^k\cBB_{\wh\Gamma}\,.
\ee
This union includes $\komoduli$ for~$\pi$ being the trivial covering of
a point to a point. Points of $\kLMS$ are called {\em \mskds}, i.e.\
the preceding definition completes the specification of the equivalence
relation stated in the introduction. Points of $\kLMS$ are thus
given by a tuple $(\wh{X},\wh{\bfz}, \wh\Gamma, \bfomega, \wh{\sigma},\tau)$
where $\bfomega = (\omega_{(-i)})_{i=0}^{L}$ is a tuple of one forms~$\omega_i$
on the subcurve corresponding to the vertices at level~$i$.
We often write just $(\wh{X},\bfomega,\tau)$ or $(\wh{X},\bfomega,\wh\Gamma)$.
The equivalence classes are given by the orbits of the action~\eqref{eq:actionLRT}
on $(\bfomega,\sigma)$.
\par
We now provide this space with a topology by exhibiting all
converging sequences. The basic idea is the conformal topology
on $\barmoduli$ where sequences converge if there is an exhaustion
of the complement of nodes and punctures and conformal maps of
the exhaustion to neighboring surfaces, see~(b) below. For multi-scale
differentials we require moreover the convergence of the differentials
as in~(c) after a rescaling, where the magnitude of rescaling is
compatible with the level structure, see~(a) and~(c). Since the
conformal topology only requires the comparison maps to be
diffeomorphisms near the nodes, which can twist arbitrarily, we
need to add~(d) to avoid constructing a non-Hausdorff space.
In the sequel we write $\wX_{\sigma_n}$ for $(\wX_n)_{\sigma_n}$ in
a sequence of welded surfaces.
\par
We say that a sequence $(\wX_n, \bfomega_n,\wh\Gamma_n)$
converges to $(\wX,\bfomega,\wh\Gamma)$, if there exist representatives of
all the equivalence classes (that we denote by the same letters),
a sequence $\ve_n \to 0$  and a sequence $\bft_n = (t_{n,i})_{i=1}^{L}
\in (\CC^*)^L$ of tuples such that:
\begin{itemize}
\item[(a)] For sufficiently large~$n$ there is an undegeneration
$\delta_n=(\delta_n^{\rm ver},\delta_n^{\rm hor})$ with
$\delta_n( \wh\Gamma)=\wh{\Gamma}_n$. 
\item[(b)] For sufficiently large~$n$ there is an almost-diffeomorphism
$g_n \colon \wX_{\bft_n \cdot \sigma} \to \wX_{\sigma_n}$ that is conformal
on the $\epsilon_n$-thick part of $(\wX, \wh\bfz)$ and that
respects the marked points, up to relabeling in~$\pi$-fibers.
\item[(c)] The restriction of $\prod_{j=1}^{i} t_{n,j}^{\ell_j} \cdot
g_n^*(\omega_{n})$ to the
$\epsilon_n$-thick part of the level~$-i$ subsurface of $(\wX, \wh\bfz)$
converges uniformly to $\omega_{(-i)}$.
\item[(d)] For any $i,j \in L(\wh\Gamma)$ with $i>j$, and any subsequence
along which $\delta^{\rm ver}_n(i)=\delta^{\rm ver}_n(j)$, we have
    \begin{equation*}
\lim_{n\to\infty} \prod_{k=|i|+1}^{|j|} |t_{n,k}|^{-\ell_k}  \= 0\,.
    \end{equation*}
\item[(e)] The almost-diffeomorphism~$g_n$ are nearly turning-number preserving,
i.e.\ for every good arc~$\gamma$ in~$\wX_\sigma$, the difference
$\rho(g_n \circ F_{\bft_n} \circ \gamma) - \rho(F_{\bft_n} \circ \gamma)$
of turning numbers converges to zero, where $F_{\bft_n}$ is the fractional
Dehn twist around the edge~$e$ by the angle $\prod_{j=1}^i t_{n,j}^{\ell_j/\kappa_e}$.
\end{itemize}
\par
This topology is exactly the topology defined in \cite{LMS} of the
compactification of the moduli spaces $\omoduli(\wh\mu)$ quotiented
by $S_{\wh{\mu}}$ and restricted
to the subspace of $k$-cyclic covers. Note that the inclusion of the covering
enhanced level graph~$\wh\Gamma$ into the datum of a \mskd implies that even
boundary points have canonically determined $k$-covers. We thus obtain
a map 
\be \label{eq:closedimm}
\kLMS \longrightarrow  \LMS[\wh{\mu}][\wh{g},\{\wh{n}\}] :=
\LMS[\wh{\mu}][\wh{g},\wh{n}] / S_{\wh{\mu}} \,.
\ee
\par
\begin{Prop}
The moduli space $\kLMS$ is a Hausdorff topological space
and its projectivization $\PP\kLMS$ is a compact Hausdorff space.
\par
The map to the partially marked stratum of abelian differentials
in~\eqref{eq:closedimm} restricted ot $\komoduli$ is a closed immersion.
\end{Prop}
\par
\begin{proof} The proof that $\kLMS$ is Hausdorff and has compact quotient
can be taken verbatim from the abelian case \cite[Theorem~9.4 and
Proposition~14.2]{kdiff}, since the topology is defined in the same way.
The second statement is local, so that we can use a marking of relative
homology and period coordinates, see e.g.\ \cite[Theorem~2.1]{kdiff}.
\end{proof}
\par
We do not enter the discussion about the local structure
of~\eqref{eq:closedimm} near the boundary, since for local computations
in Section~\ref{sec:good} we use the model domain we now define
and the resulting perturbed period coordinates.

\subsection{Model differentials and modifying differentials}
\label{sec:Model}

In order to provide $\kLMS$ with a complex structure we use a
local model space that automatically has a complex structure
(as a finite cover of a product of spaces of non-zero $k$-differentials).
The degeneration of differentials on lower components is
emulated in the model space by vanishing of auxiliary parameters~$t_i$.
\par
A first attempt would be to consider that action~\eqref{eq:actionLRT} of
the level rotation torus $\LRT$, that  makes $\tkdpm$ into a principal
$(\CC^*)^{L}$-bundle over the quotient space. However $\LRT$ is in general
not naturally isomorphic to $(\CC^*)^{L}$ and so there is no natural
associated $\CC^{L}$-bundle that could serve as local compactification.
But the simple  level rotation torus $\sLRT$ has such an isomorphism,
as remarked along with~\eqref{eq:sLRTpara}.
\par
The successful strategy is to construct a $K_{\wh{\Gamma}}:= \Tw/\sTw$-cover of $\tkdpm$
on which now $\sLRT$ acts as a lift of the action ~\eqref{eq:actionLRT}.
In \cite[Section~5]{LMS} such a space is constructed using an additional
Teichm\"uller marking, and here we use the subspace where the flat surfaces
admit a~$\tau$-action. On this space even the universal cover of $\LRT$ acts,
and one can then quotient by $\sTw$ and forget the marking of the rest of the
surface to obtain the requested \emph{(uncompactified) simple model domain}
$\tkdpms$.
\par
The action of $\sLRT$ on $\tkdpms$ make this space a principal
$(\CC^*)^{L}$-bundle over the quotient space, which we call with hindsight 
the 'simple' version of the boundary stratum $\Omega^k\cBB^s_{\wh\Gamma}
= \tkdpms / \LRT^s$. We define the {\em simple model domain} to be the
associated $\CC^{L}$-bundle, 
\be \label{eq:assCbundle}
\overline{\tkdpms} \= (\CC^{L} \times \tkdpms) / \sim,
\qquad (\bft, \bfeta) \sim (\rho \cdot \bft, \rho^{-1} \cdot \bfeta) \quad
\text{for $\rho \in \sLRT$}
\ee
over $\Omega^k\cBB^s_{\wh\Gamma}$. The smoothness
of strata of $k$-differentials and the smoothness of the associated $\CC^{L}$-bundle
directly implies:
\par
\begin{Prop}
The compactified simple model domain $\overline{\tkdpms}$ is smooth
with normal crossing boundary divisor given by the divisors $D_i = \{t_i = 0\}$
with $t_i$ the local coordinates on $\sLRT$ as in~\eqref{eq:sLRTpara}.
\end{Prop}
\par
Note that the boundary of the compactified simple model domain comes
with a natural {\em stratification} given by the subset of
$\{1,\ldots,L\}$ of the~$t_i$ that are zero.
\par
Our next goal is to exhibit the universal curve and a universal family
of differentials over $\overline{\tkdpms}$. After adding modifying 
differentials and performing a plumbing construction this gives a family
of multi-scale differentials in the proof of the main theorem below.
\par
A change of notation seems adequate to illustrate the process: So far
we denoted a twisted $k$-differentials  by $(\wh{X}, \bfomega, \tau)$, since
they arose from limits of canonical covers of $k$-differentials. Subsequently
we will however denote curves with differentials over the simple model domain
$\overline{\tkdpms}$ by~$(\wh{Y},\bfeta)$ (as we did already in~\eqref{eq:assCbundle})
and call them \emph{model differentials}, even though the space $\overline{\tkdpms}$
is made from a space of twisted differentials
by adding prong data and compactification. We will then call $(X,\bfomega)$
the curves with differentials obtained from plumbing. These will yield
multi-scale $k$-differentials. To memorize: the model differentials~$\bfeta$
are always non-zero, the curves are equisingular and the boundary is specified by
vanishing of the auxiliary parameters~$t_i$. The differentials $\bfomega$ will be
constructed below on degenerating curves, they will vanish on lower level components
and the boundary is given by the appearance of nodes.
\par
The space $\tkdnoab$, being the subspace defined by global residue
conditions  in a product of moduli spaces, obviously comes with a universal family
of curves and $k$-differentials that we can pull back to $\tkdpm$.
Since the level rotation torus only acts on differentials and
prong-matchings, not on the curve, the universal curve descends
to a family of curves  $f\colon \wh\cYY  \to \overline{\tkdpms}$.
\par
Next we turn to the differentials. It will be convenient to fix the scale
of the $\sLRT$-orbits of the differentials in~\eqref{eq:assCbundle} rather
than working with equivalence classes. That
is, over small enough open sets~$W \subset \overline{\tkdpms}$ we may
assume that we work with a fixed collection $\bfeta = (\eta_{(-i)})_{i=0}^{L}$ of
families of differentials representing in each fiber the equivalence class
of the corresponding point in $\Omega^k\cBB^s_{\wh\Gamma}$. Consequently the
collection of functions $\bft = (t_i)$ is part of a coordinate system on~$W$.
\par
\medskip
The modifying differentials we now define will be used for plumbing and also
for perturbed period coordinates on charts of $\overline{\tkdpms}$.
In the sequel we check that the construction of \cite[Section~11]{LMS}
works in the $k$-equivariant setup. We define
\begin{equation} \label{eq:omTeta}
\bft\ast \bfeta \= \left(\prodt\cdot \eta_{(-i)} \right)_{i=0}^L \=
\left(t_{1}^{\ell_{1}}\dots t_{i}^{\ell_{i}} \cdot \eta_{(-i)} \right)_{i=0}^L \,,
\end{equation}
for $\bft=(t_{1},\dots, t_{L})\in (\CC^*)^L$ and use  by definition the
trivial rescaling $t_0=1$ on top level unless specified otherwise.
\par
\begin{Defi}\label{df:modif}
An \emph{equivariant family of modifying differentials} on the universal family 
$f\colon \wh\cYY\to W$ restricted to~$W$
equipped with the universal differential $\bft \ast \bfeta$  is a
family of meromorphic differentials $\bfxi = (\xi_{(-i)})_{i=0}^{L}$ on
$f\colon \wh\cYY\to W$, such that
\begin{enumerate}
\item[(i)] the equivariance $\tau^* \bfxi = \zeta \cdot \bfxi$ holds, 
\item[(ii)] the differentials $\xi_{(-i)}$ are holomorphic, except for possible
simple poles along both horizontal and vertical nodal sections, and except
for marked poles,
\item[(iii)] the component $\xi_{(-L)}$ vanishes identically and moreover
$\xi_{(-i)}$ is divisible by $t_{\lceil i+1\rceil}$ for each $i=1,\ldots,L-1$, and
\item[(iv)] the sum $\bft\ast \bfeta  + \bfxi$ has opposite residues at every
 node.
\end{enumerate}
\end{Defi}
\begin{Prop}   \label{prop:constrModif}
The universal family $f \colon \wh\cYY \to W$ equipped with the universal
differential $\bft \ast \bfeta$ admits an equivariant family of modifying
differentials.
\end{Prop}
\par
\begin{proof} The proof of \cite[Proposition~11.3]{LMS} works in
the situation where the edges of $\wh\Gamma$ are images of the pinched
multicurve~$\Lambda$ via a family of markings $\Sigma \to \wY$
by a reference surface~$\Sigma$. Choosing~$W$ contractible, we may assume
that we have such a marking here as well.
\par
The proof in loc.\ cit.\  starts by taking the subspace 
$V = \langle \lambda \in \Lambda \rangle_\QQ$ and the subspace~$V_P$ spanned by
the loops around the marked poles inside $H_1(\wY \smallsetminus \wh{P},\QQ)$.
We define $V_N=V+V_P$. The proof proceeds by  searching for a complementary
subspace~$V_C$ (i.e. with $V_N \cap V_C = 0$) 
such that the projection $p(V')$ to $H_1(\wY, \QQ)$ is a Lagrangian subspace,
where $V' = V_N + V_C$. The proof then constructs~$\bfxi = (\xi_{(i)})$
for each $w \in W$ from assignments  $\rho_i: V_i \to \CC$ determined
by the periods of the fiber~$\bfeta_w$ on subspaces $V_i$ of~$V+V_P$
generated by multicurves associated to edges whose lower level is below~$i$.
Relevant here is that those~$\bfxi$ satisfy all properties of Definition~\ref{df:modif}
except possibly the equivariance in~(i). Moreover, $\bfxi$ depends uniquely
on an extension~$\rho'_i$ of $\rho_i$, that we may chose to be zero on $V_C$.
\par
If we can find a subspace $V_C$ which is $\tau$-invariant, then the
extended residue assignment~$\rho'_i$ is $\tau$-equivariant (with $\tau$ acting
by multiplication by~$\zeta$ on the range) and thus $\bfxi$ satisfies~(i).
To find such a~$V_C$, we enlarge~$V_C$ and thus $V'=V_C+V_N$ step by step,
staying $\tau$-invariant at each step, until $p(V')$ is a Lagrangian subspace.
If at some step $V'$ is $\tau$-invariant, but $p(V')$ is strictly contained in
a Lagrangian subspace, we may find an element $\gamma$ that pairs trivially
with $p(V')$. But then $\tau^i(\gamma)$ also pairs trivially with $p(V')$ for
all~$i$ and we add to $V_C$ the span of all $\tau^i(\gamma)$.
\end{proof}

\subsection{The perturbed period map}
\label{sec:PPer}

Periods give local coordinates on $\tkd$ and thus on $\tkdpm$.
Together with the tuple of degeneration parameters~$\bft$ and deprived
of one coordinate per level to fix the scale of projectivization they
give local coordinates of $\overline{\tkdpms}$. We introduce some
perturbation of these
coordinates here and show that this still gives local coordinates.
The reason for this procedure is that the perturbed period coordinates
can still be used after plumbing, see Section~\ref{sec:ProofThm11}. Together
with horizontal plumbing parameters it will provide coordinates on an
orbifold chart of~$\kLMS$. Except for the use of appropriate eigenspaces
this is exactly \cite[Section~11]{LMS}.
\par
We use for the rest of the discussion a small enough~$W$ as above and
a fixed modifying differential provided by Proposition~\ref{prop:constrModif}
after choosing a curve system, as in its proof.
Near the marked point $e^+$ corresponding to the upper end (say on
level~$i = i(e^+)$) of each of the vertical nodes, choose an auxiliary
section $s_e^+: W \to \wh\cYY$ such that
\be \label{eq:constdistance1}
\int_{e^+}^{s_e^+(w)} \bfeta_{(i)} \= {\rm const}\,,
\ee
where the constant is sufficiently small, depending on~$W$, and
constrained by the plumbing construction later. Near each zero marked~$z_j$
of~$\bfeta$ (say on level $i=i(z_j)$) choose an auxiliary section
$s_j: W \to \wh\cYY$ that coincides with the barycenter of the zeros
of $\eta_{(i)} + t_{\lceil -i\rceil}^{-1} \cdot\xi_{(i)}$ that result from the deformation
of $z_j$. We let $\gamma_{ij}$ for $i=0,\ldots,-L$ and $j=1,\ldots,
\dim E_{(i)}^{\rm grc}$ be a basis of the subspaces $E_{(i)}^{\rm grc}$ of
homology. Since the contribution of each level to the twisted differentials
compatible with a level graph is positive-dimensional (by the rescaling
of the differential), the definition of periods coordinates along the boundary
in Proposition~\ref{prop:dimper} implies that for each~$i$ there exists
some $j$ such that $\int_{\gamma_{i,j}} \eta_{(i)} \neq 0$.  We use this
to fix the scale of the projectivization and assume that the periods
for $j=1$ are normalized on each level, i.e. $\int_{\gamma_{i,1}} \eta_{(i)} =1$.
\par
The $i$-th level component of the perturbed period map is now given by 
\par\be \label{eq:defPPer}
\PPer_i\colon \begin{cases} \begin{array}{ccl} W &  \to &
\CC^{\dim E_{(i)}^{\rm grc} -1 + \delta_{i,0}}\,,  \\
[(\wY,\bfeta,\bft)] &\mapsto & \left(\int_{\gamma_{i,j}}  \eta_{(i)}
+  t_{\lceil -i\rceil}^{-1} \cdot\xi_{(i)}\right)_{j=2-\delta_{i,0}}^{\dim E_{(i)}^{\rm grc}}\,,
\end{array}
\end{cases}
\ee
where $\delta_{i,0}$ is Kronecker's delta and where the integrals are to
be interpreted starting and ending at the nearby points determined by
the sections~$s_e^+$ and $s_j$ rather than the true zeros of~$\bfeta$.
The reason for this technical step is that those nearby points are
still present after the surfaces has been plumbed ('Step~2' below).
Recall that we defined $D_i = \{t_i = 0\}$.
\par
\begin{Prop} \label{prop:pertper}
The perturbed period map 
$$\PPer_{\wh{\Gamma}}^{\rm MD}\colon W \to \CC^{L} \times \prod_{i = 0}^{-L}
\CC^{\dim E_{(i)}^{\rm grc} -1 + \delta_{i,0}}\,,
\quad [(\wY,\bfeta,\bft)] \mapsto  \Bigl(\bft\,;\,\coprod_{i = 0}^{-L}
\PPer_i(\wY,\bfeta,\bft)\Bigr)$$ 
is open and locally injective on a neighborhood of the most degenerate
boundary stratum $W_{\wh{\Gamma}} = \cap_{i=1}^{L} D_i$  in the
compactified model domain $\overline{\tkdpms}$.
\end{Prop}
\par
We will write $(\bft,\bfw)= \PPer_{\wh{\Gamma}}^{\rm MD} (\wY,\bfeta,\bft)$.
\par
\begin{proof}
As in \cite[Proposition~11.6]{LMS} it suffices to show that the
derivative is surjective, by dimension comparison. For the tangent
directions to the boundary this follows from Proposition~\ref{prop:dimper}
(and the fact that we have projectivized the lower levels). For the
transverse direction this follows since the~$t_i$ are the local coordinates
of the $\CC^{L}$-bundles used to construct the compactifications.
\end{proof}
\par
The reader should keep in mind, that in the model domain with its equisingular
family of curves horizontal nodes are untouched. They enter in
Proposition~\ref{prop:PPerIntro} only after plumbing horizontal nodes,
see below.

\subsection{The complex structure and the proof of Theorem~\ref{thm:kLMS}}
\label{sec:ProofThm11}

The outline of the proof of Theorem~\ref{thm:kLMS} consists of the
following steps.
\begin{itemize}
\item[1)] Construct locally covers~$U^s \to U$ for small open sets
$U \subset \kLMS$ that will be used as orbifold charts.
\item[2)] Perform a plumbing construction on the pullback
of the universal family $f\colon \wh\cYY  \to \overline{\tkdpms}$ to
small open sets~$W$ and via the second projection to $W \times \Delta_\ve^h$ 
in order to obtain a family $\cXX \to W \times \Delta_\ve^h$ together with
a family of differentials. 
\item[3)] Use the moduli properties of the strata of $\kLMS$ to
obtain moduli maps $\Omega{\rm Pl}: W \times \Delta_\ve^h \to U^s$
 for an appropriately chosen target set~$U^s$, defined stratum by stratum.
\item[4)] Show that $\Omega{\rm Pl}$ is a homeomorphism near
a central point $P \times (0,\ldots,0) \in W \times \Delta_\ve^h$
and thus provide charts there.
\end{itemize}
\par
The charts constructed in this way depend on many choices, in the construction
of the modifying differential and the parameters for plumbing. However,
the induced complex structures fit together and that's all we need since
$\kLMS$ already exists as a topological space. We provide the
details for  Step~2), since there the $\tau$-equivariance
needs to be respected and since we need this in the next section.
The technical Step~1),  Step~3) and Step~4) proceed exactly as in \cite{LMS}.
\par
\medskip
\paragraph{\bf Step~1} In order to provide $\kLMS$ with a complex structure
we consider the neighborhood~$U$ of a point $(\wh{X},\wh{\bfz}, \bfomega,
\tau, \wh\Gamma)$ that we may assume to be at the boundary, say for level
graph~$\wh\Gamma$.
(The following description assumes that $(\wh{X},\bfomega,\tau)$ has no automorphisms.
In general we should start from an orbifold chart, and add the extra orbifold
structure described below.)
The compactified simple model domain is a $K_{\wh{\Gamma}} =\Tw/\sTw$-cover of 
the (in general) singular space that we would get by compactifying
the $\LRT$-quotient of $\tkdpm$. Consequently, we have to pass
locally near $(X,\bfomega, \wh\Gamma)$ to a~$K$-cover of~$U$. We
define this cover~$U^s$ as follows. Define the auxiliary {\em simple
boundary stratum} to be $\Omega^k\cBB^s_{\wh\Gamma} = \tkdpms / \sLRT$.
As a set
\bes
U^s = \Bigl\{  (X',\bfomega', \wh\Gamma') \in 
\bigcup_{\wh\Gamma' \rightsquigarrow \wh\Gamma}\Omega^k\cBB^s_{\wh\Gamma'}
\,\,\colon \,\,\varphi((X',\bfomega', \wh\Gamma')) \in U
\Bigr\} \,,
\ees
where $\varphi$ is induced by the quotient maps $\Omega^k\cBB^s_{\wh\Gamma'}
\to \Omega^k\cBB_{\wh\Gamma'}$. We provide $U^s$ with a topology where
convergence is formally given exactly by the same conditions as for~$\kLMS$
in Section~\ref{sec:TopoComp}, but where now the 'existence of
representatives of the equivalence classes' is up to the
torus~$\sLRT$ rather than the quotient torus~$\LRT$. 
\par
\medskip
\paragraph{\bf Step~2} To start the plumbing construction  we first 
define the {\em plumbing fixture} for each vertical edge $e \in E(\wh\Gamma)$
to be the degenerating family of annuli
\be \label{eq:plumbing_fixture}
 \VV_e \= \Bigl\{(\bfw, \bft, u, v) \in W \times \Delta_{\delta}^2 : 
uv = \prod_{i =-e^++1}^{-e^-} t_i^{m_{e,i}}\Bigr\},
\ee
 that only depends on the $\bft$-part of the perturbed period coordinates
$(\bft, \bfw)$ of~$W$. Recall that by definition \eqref{eq:defai} we have set $m_{e,i} \= \ell_i/\kappa_e$. We equip $\VV_e$ with the family of differentials
\be \label{eq:reloneform}
\Omega_e \= \left(t_{\lceil -e^+\rceil}\cdot u^{\kappa_e -1} - \frac{r'_e}{u}\right)
du \= \left(- t_{\lceil -e^-\rceil}\cdot v^{-\kappa_e -1} + \frac{r'_e}{v}\right)dv\,,
\ee
where we recall that $\prodt  \= t_i^{\ell_i} \dots \,t_1^{\ell_1}$
and where $r'_e = r'_e(\bfw,\bft)$
are the residues of the universal family over model domain. Inside the plumbing
fixture we define the gluing annuli $\cAA_e^\pm$ by $\delta/R <|u| < \delta$ 
and $\delta/R <|v| < \delta$ respectively. The sizes~$\delta$, $R$ and the
size of the neighborhood~$W$ will be determined in terms of the geometry
of the universal family, to ensure for example that plumbing
annuli are not overlapping.
\par
Suppose we only have vertical nodes. The plumbing construction proceeds bottom
up. Near each of the nodes of bottom level we put the family of differentials
$\eta_{(-L)}$ in standard form $(v^{-\kappa_e -1} + \tfrac{r_e}{v})dv$ so that 
after rescaling with $t_{\lceil -e^-\rceil}$ it can be glued to~$\Omega_e$ for $r_e' = t_{\lceil -e^-\rceil} r_e$.
That such a normal form exists in families is the content of
\cite[Theorem~3.3]{LMS}. The functions $r_e'$ determine the modifying
differential~$\xi_{(-L+2)}$ as the proof of Proposition~\ref{prop:constrModif} 
shows, see \cite[Corollary~11.4]{LMS}. We will thus put 
$t_{\lceil -e^+\rceil} \eta_{(-L)} + \xi_{(-L+1)}$ in standard form near~$e^+$ using the
normal form on the deformation of an annulus (\cite[Theorem~4.5]{strata} or
\cite[Theorem~12.2]{LMS}) and this glues with the form~\eqref{eq:reloneform}
on the upper end of the annulus. Iterating the procedure allows to plumb the
collection of families of one-forms
\be \label{eq:omTetaxi}
\bft\ast \bfeta + \bfxi \= \left(\prodt\cdot \eta_{(-i)} + \xi_{(-i)}
\right)_{i=0}^{L} 
\ee
on the equisingular family of curves $\cYY \to W$ to a family of
one-forms~$\bfomega$ on a degenerating family of curves $\cXX \to W$.
The zeros of higher order of~$\bfeta$ may have split up in~$\bfomega$
when adding~$\bfxi$. A local surgery merges them back to the barycenter
\cite[Lemma~4.7]{strata}.
\par
In the preceding construction we have neglected so far that the choice
of the normal form is unique only up multiplication by a $\kappa_e$-th root 
of unity. The prong-matching that is part of the datum of the universal family
over the model domain determines this choice. Many more details, using reference
sections to make the construction rigorous, are given in \cite[Section~12]{LMS}.
\par 
The whole construction can obviously performed $\tau$-equivariantly, 
since the modifying differential is $\tau$-equivariant and since the
sizes of the neighborhoods and plumbing annuli are determined by the
rates of degeneration of $\bft \ast \bfeta + \bfxi$, i.e.\ by $\tau$-equivariant 
data.
\par
Finally, we investigate horizontal nodes of $\wh\Gamma$, that come in
$\tau$-orbits of length~$k$ and that we thus label as $e_1^{(a)},\ldots,
e_{|E^h|}^{(a)}$ for $0 \leq a <k$. We parameterize the plumbing by additional
plumbing parameters $\bfx = (x_1,\ldots,x_{|E^h|}) \in  \Delta_\epsilon^{|E^h|}$ and
define the  {\em (horizontal) plumbing fixture} to be
\begin{equation}\label{eq:horplumbfix}
  \WW_j \= \left\{ (\bfw, \bft, \bfx, u, v) \in W \times
    \Delta_\epsilon^{|E^h|} \times \Delta_\delta^2 : uv=x_j\right\},
\end{equation}
independently of the upper label~$a$ of $e_j^{(a)}$, equipped with the
family of holomorphic one-forms
\begin{equation}\label{eq:horplumbform}
  \Omega_j \= -r_{e_j}'(\bfw, \bft) du/u \= {r_{e_j}'(\bfw, \bft)}dv/v\,,
\end{equation}
where $\pm r_{e_j}'(\bfw, \bft)$ is the residue of  $\bft \ast \bfeta + \bfxi$ at the 
$j$-th horizontal node. Here the gluing happens along annuli
$\cBB_j^\pm$ by $\delta/R <|u| < \delta$ and $\delta/R <|v| < \delta$.
\par
\medskip
\paragraph{\bf Step~3} The existence of moduli maps
on each stratum of the simple model domain to~$\kLMS$ is immediate from
the construction  of $\kLMS$ as union of  strata
$\Omega^k\cBB_{\wh\Gamma} = \tkdpm / \LRT$
and the property of $\tkdpm$ as moduli space of $k$-differentials.
We let~$U$ be the range of the union of these maps. The map factors
through~$U^s$ since both this space and the simple model domain are defined 
as $\sLRT$-equivalence classes. \cite[Section~12.5]{LMS} provides more
details.
\par
\medskip
\paragraph{\bf Step~4} To show that the resulting map $\Omega{\rm Pl} \colon 
W \times \Delta^{|E^h|} \to U^s$ is continuous we have to invoke the definition 
of the topology on~$U^s$ to show that the images of a converging sequence
converges. This entails exhibiting the almost-diffeomorphisms~$g_n$ with
the properties (a)--(e). These~$g_n$ are construct level by level, bottom up, 
using conformal identifications of flat surfaces with the same periods
(\cite[Theorem~2.7]{LMS}), a $C^1$-quasi-conformal extension of these maps
across the plumbing cylinder and the equivalence of the conformal and
$C^1$-quasi-conformal topology on strata of abelian differential
(\cite[Section~2]{LMS}).
\par
To show that  $\Omega{\rm Pl}$ is a homeomorphism we need to show that this
map is open and locally injective. Openness amounts to showing that for any
converging sequence in~$U^s$, say converging to $(\wh{X},\wh{\bfz},\omega,\tau,
\wh\Gamma)$, we can eventually undo the plumbing construction and find
$\Omega{\rm Pl}$-preimages in the model domain $\overline{\tkdpms}$.
These preimages are again found level by level, the scales~$t_i$ of the
model differentials being determined by the scales~$s_i$ in the definition
of convergence in~$U^s$.  Local injectivity amounts to checking uniqueness
of the previous unplumbing steps using perturbed period coordinates
and can be checked applying \cite{LMS}. See Section~12.5-12.7 in loc.\ cit.\
for details on these steps.
\par
\medskip
The action of $\CC^*$ on the $k$-th root~$\omega$ defines an action on the
space $\komoduli$ that is equivariant via $\lambda \mapsto \lambda^k$ with
a $\CC^*$-action on $\komoduli$. The quotients of both actions is the
same space $\PP \kLMS$. We encourage the reader to revisit all the 
steps to check that the first action extends equivariantly to all auxiliary
spaces, multiplying simultaneously all forms at all levels. The resulting
quotient of $\kLMS$ by $\CC^*$ is the compactification $\PP\kLMS$ claimed
in iii) of Theorem~\ref{thm:kLMS}. 
\par
\medskip
The {\em proof of Proposition~\ref{prop:PPerIntro}} is contained in these
statements, since Proposition~\ref{prop:pertper} together with the
disc coordinates~$x_j$ used in~\eqref{eq:horplumbfix}  gives local coordinates
on $\overline{\tkdpms} \times \Delta^{|E^h|}$. Consequently, the
perturbed period coordinates are given by
\ba \label{eq:PPerexpl}
\PPer \colon \quad 
U^s & \xrightarrow{\Omega{\rm Pl}^{-1}} W\times \Delta^{|E^h|}
&&\longrightarrow \CC^h \times \CC^{L+1} \times \prod_{i=0}^{-L} \CC^{\dim E_{(i)}^{\grc} -1 }\,\\
\bigl[\wX,\bfomega\bigr] & \xmapsto{\Omega{\rm Pl}^{-1}} [(\wY,\bfeta,\bft,\bfx)] &&\mapsto  \Bigl(\bfx\,;\bft\,;\,\coprod_{i = 0}^{-L}
\PPer_i(\wY,\bfeta,\bft)\Bigr)
\ea
on open orbifold charts~$U^s$ of $\kLMS$, using the inverse of the
homeomorphism $\Omega{\rm Pl}$ constructed in Step~$3$ and~$4$. 
\par

%% file: sec_good_new.tex
\section{The area form is good enough}
\label{sec:good}

Here we prove our main Theorem~\ref{thm:areagood}. We place ourselves
in the setting of the theorem and recall that now $m_i >-k$ and thus
the sets~$P$ and $\wh{P}$ as defined in Section~\ref{sec:period} are empty.
The first step is to  determine where the metric tends to infinity and
then to give a convenient expression of the metric. Arguing inductively
on~$k$, we may also suppose that we are dealing with primitive
$k$-differentials, i.e.\ that the canonical $k$-cover is connected.
\par
We start with the definition of the corresponding hermitian form. For a
symplectic basis $\alpha_1,\ldots,\alpha_{\wh{g}}, \beta_1,\ldots,\beta_{\wh{g}}$
of the absolute homology $H_1(\wh{\Sigma},\ZZ)$ and for $\omega,\eta
\in H^1(\wX,\CC)$ we define hermitian form
\be
\langle \omega,\eta \rangle \= \frac{i}{2}
\sum_{i=1}^{\wh{g}}  \bigl(\omega(\alpha_i)\,
\overline{\eta(\beta_i)} - \omega(\beta_i)\,\overline{\eta(\alpha_i)}\bigr)
\ee
with the abbreviations $\omega(\alpha_i)  = \int_{\alpha_i} \omega$ etc.
By Riemann's bilinear relations we can rewrite the metric defined
in~\eqref{eq:defh} using the above hermitian form as 
\bes
h(X,q)^{1/k} \= \langle \omega,\omega \rangle
\= \frac{i}{2} \sum_i (a_i \overline{b_i} - b_i \overline{a_i})\,,
\ees
where we introduce another abbreviation $a_i = \omega(\alpha_i)$
and $b_i = \omega(\beta_i)$, to be used if~$\omega$ is the only one-form
that appears. We recall from~\eqref{eq:omTeta} the notation
\[\tui=t_1^{\ell_1}\cdots t_i^{\ell_i}\]
for any $i=0,1,\dots,L$ (in particular $t_{\lceil 0\rceil}=1$).
\par
\begin{Lemma} \label{le:smoothvert}
  On a neighborhood $U$ of a boundary point whose corresponding level
graph~$\wh\Gamma$ has only vertical edges and $L+1$ levels,
the metric~$h$ extends to a function of the form
\begin{equation}\label{eq:metricvertical}
  h(X,q)^{1/k}\= \sum_{i=0}^{L} \left|t_{\lceil i\rceil}\right|^2
\Bigl(h^{\text{tck}}_{(-i)}-\sum_{p=1}^{i}R^{\text{ver}}_{(-i),p} \log\left|t_p\right|
\Bigr)\,,
\end{equation}
where $h^{\text{tck}}_{(-i)}$ is a smooth positive function bounded away from zero and $R^{\text{ver}}_{(-i),p}$ is a  smooth non-negative function.
\end{Lemma}
\par
\begin{proof}   A neighborhood of the point~$(\wX,\bfomega)$ is also
a neighborhood~$U$ of that point in the model domain. There, $\bfomega = \bfeta$
and it is interpreted as a collection $\eta_{(-i)}$ of non-zero
differential forms on the subsurface $\wX_{(-i)}$ on the $-i$-th level. The
neighborhood of~$(\wX,\bfomega)$ consists of the stable differentials obtained
via the plumbing construction applied to the differential forms
$(\prod_{j=1}^{i} \,t_j^{\ell_j}) \eta_{(-i)} + \xi_{(-i)}$ on the universal
family over model domain restricted to the small neighborhood. Here $\bft =
(t_i)_{i=1}^{L}$
is the collection of 'opening-up' parameters in the polydisc and the positive
integers~$\ell_j$  are determined by the enhanced level graph~$\Gamma$
via~\eqref{eq:defai}. The central fiber of this family agrees
with~$\bfomega$ by construction. 
\par
We compute the local expression of the metric in the neighborhood~$U$.
We decompose a plumbed surface $(\wX_u,\omega_u)$ over $u\in U$ in the
following way. Let $\wh\cXX_{(-i),u}$ be  the fiber  over $u \in U$ of
the $-i$-th level subsurface over the model domain. For any edge~$e$ connected
to level~$-i$, consider the complement of the interior of the gluing
annuli~$\cAA_e^{+/-}$ (defined in Step~$2$ above), so that we remove
neighborhoods of the points where $\eta_{(i)}$ has zeroes or poles
corresponding to edges of $\wh\Gamma$.  We thus decompose the plumbed
surface  as 
\[\wX_u \=\bigsqcup_{i=0}^{L} \wh\cXX_{(-i),u}^\circ\quad  \sqcup \quad \bigsqcup_{e\in E(\wh\Gamma)} \VV_e \,,\]
where $\VV_e$ are the plumbing fixtures at the edges defined in \eqref{eq:plumbing_fixture}. We have
\be \label{eq:arearough}
\area_{\wX_u} (\omega_u) \=\sum_{i=0}^{L} \area_{\wh\cXX_{(-i),u}^\circ}(\prodt 
\eta_{(-i),u} + \xi_{(-i),u})\ +\ \sum_{e\in E(\wh\Gamma)} \area_{\VV_e}\Big(\Omega_e\Big)
\,,\ee
where the differential form $\Omega_e$ was defined in~\eqref{eq:reloneform}.
The first summands above give a smooth function in $U$. Moreover, since the
components $\xi_{(-i),u}$ of the modification differentials are divisible
by $t_{\lceil i+1\rceil}$ (see Definition~\ref{df:modif}~iii)), we can write
\[ \area_{\wh\cXX_{(-i),u}^\circ}(\prodt 
\eta_{(-i),u} + \xi_{(-i),u}) \=\left|t_{\lceil i\rceil}\right|^2h^{\text{tck}}_{1,(-i)}\]
where $h^{\text{tck}}_{1,(-i)}$ is a smooth function bounded away from zero on $U$.
\par
It remains to compute the second summands in the expression~\eqref{eq:arearough}.
For each edge $e$ of the graph $\wh\Gamma$, let
$t_{\Delta(e)}:=\prod_{p =-e^+ +1}^{-e^-} t_p^{m_{e,p}}$, where recall from \eqref{eq:defai} that
$m_{e,p}=\ell_p/\kappa_e$. Hence we obtain, using $r_e = r_e'/\prodtbot$
as defined after~\eqref{eq:reloneform}, that 
\ba \label{eq:areaVV}
  \area_{\VV_e}\Big(\Omega_e\Big)
&=\frac{i}{2}\int \limits_{\frac1\delta {|t_{\Delta(e)}|}
\leq |z|\leq \delta} \left|\prodttop\cdot z^{\kappa_e -1}
- \frac{\prodtbot r_e(u)}{z}\right|^2 dzd\bar{z}\\
&=2\pi  \int_{\frac1\delta {|t_{\Delta(e)}|}}^\delta
\left(\left|\prodttop\right|^2 r^{2(\kappa_e -1)}+\left|\prodtbot\right|^2 |r_e(u)|^2 \frac{1}{r^2})\right)rdr\\
&\phantom{\=\,} -\left|\prodttop\prodtbot\right| |r_e(u)|^2\int_{0}^{2\pi}
\int_{\frac1\delta {|t_{\Delta(e)}|}}^\delta r^{\kappa_e-2}\Big(e^{i\kappa_e  \theta}+e^{-i\kappa_e  \theta}\Big)r dr d\theta\\
&=\left|\prodttop\right|^2 h^v_e(u)-2\pi |r_e(u)|^2 \left|\prodtbot\right|^2 \sum_{p=-e^+ +1}^{-e^-}m_{e,p}\log\left|t_p\right|
\ea
where $h^v_e(u)$ is a smooth function on $U$. Here the third line integrates
to zero, the lower integration limit of the $dr/r$-integral
in the second line gives the last term and all the rest has been subsumed
in $h^v_e(u)$.
\par
We now rearrange the sum over all the vertical plumbing fixtures contributions
according to the level of the top end of the edge. Moreover we group together
the sum of the areas of the thick part and the non-residue part of the plumbing
fixtures to get an expression
\[h_{(-i)}^{\text{tck}}\,:= \, h^{\text{tck}}_{1,(-i)}+\sum_{e \colon e^+ = -i} h^v_e(u)\,.\]
This is a smooth function bounded away from zero. We also group the residue
terms, the last terms in~\eqref{eq:areaVV}, but now according to the
bottom end of the edge. We set for all $0 < p \leq i$: 
\[ R^{\text{ver}}_{(-i),p} \= -\sum_{e \colon -e^- = \,i, \, p < -e^-}
2\pi |r_e(u)|^2 m_{e,p} \,. \]
This gives~\eqref{eq:metricvertical}.
\end{proof}
\par
Now suppose we work in a neighborhood~$U$ of a  boundary point with only horizontal
nodes. Assume there are $E^h = E^h_{(0)}$ horizontal nodes and
let $x_j$ for $j=1,\dots,E^h$ be their opening-up parameters.
\par
\begin{Lemma}\label{le:smoothhor}
On the neighborhood~$U$ the metric~$h$ has the form
\begin{equation}\label{eq:metrichor}
h(X,q)^{1/k} \= h^{\text{tck}}_{(0)} \,-\,
\sum_{j=1}^{E_{(0)}^h} R^{\text{hor}}_{(0),j} \log(|x_{j}|^2)\,,
\end{equation}
where $h^{\text{tck}}_{(0)}$ and $R^{\text{hor}}_{(0),j}$ are smooth functions independent
of the $x_j$ parameters, both bounded above and away from zero.
\end{Lemma}
\par
\begin{proof} The total space of the line bundle $\cOO(-1)$ defined in the
introduction as the $\varphi$-pullback of $\cOO(-1)$ from the incidence variety
compactification, is nothing but the total space of the projection
$\kLMS \to \PP\kLMS$. Our goal is thus to find an expression for the
area of a point in~$\kLMS$  near a boundary point
$(\wX,\bfomega,\wh\Gamma) \in \partial \PP\kLMS$.
\par
For notational simplicity we consider first the case {\bf that $X$ has only
one horizontal node that we moreover suppose to be non-separating}. Consequently,
$\wX$ has~$k$ nodes. We pick a convenient basis of $H_1(\wh{\Sigma},\ZZ)$ on
a smooth model $\wh{\Sigma}$ (connected by our standing primitivity
assumption) that is pinched to~$\wh{X}$. The $k$ pinched curves $\alpha_i
\in \wh\Sigma$ are linearly independent and form a $\tau$-orbit in homology. 
Next, we take the symplectic dual curves $\beta_{i}$ with the intersection
pairing $\langle \alpha_{i}, \beta_{j}\rangle = \delta_{ij}$.
Note that $\beta_{i}$ is well-defined in a neighborhood of
$(\wX,\bfomega,\wh\Gamma)$ (only) up to adding an integer multiple
of~$\alpha_{i}$.  We arbitrarily complement
these elements  by $\alpha_i,\beta_i \in H_1(\wh{\Sigma},\ZZ)$
for $i=k+1,\ldots,\wh{g}$ to a symplectic basis.
\par
In the current case the multi-scale differential case
$\bfomega = (\omega_0)=(\eta_0)$
consists of a single one-form. Recall from Step~2 in Section~\ref{sec:ProofThm11}
that points in a neighborhood of $(\wX,\bfomega,\wh\Gamma)$ are obtained from
surfaces $(\wX',\eta') \in \omoduli[\wh{g}-k,\wh{n}+2k](\wh\mu, (-1)^{2k})$
that admit an action by $\langle \tau \rangle \cong \ZZ/k$, 
by gluing in $k$~times each of the plumbing fixtures~$\WW$ in a
$\tau$-equivariant way, parameterized by a coordinate
$\bfx = (x) \in \Delta$ as in Step~2 above. By Proposition~\ref{prop:PPerIntro}
and explicitly~\eqref{eq:PPerexpl} the coordinates near the boundary point
are~$\bfx$ and the periods in the $\zeta_k$-eigenspace of~$\eta'$. We denote
by $\omega'$ the  differential obtained from $\eta'$ after the plumbing construction.
Notice that $\omega'$ is a holomorphic differential on the plumbed surfaces having
all plumbing parameters $x_i$ different from zero.  Our aim is to rewrite the area
form, which is defined using $\omega'$ periods in the interior, in terms of the
perturbed period coordinates, i.e. $\bfx$ and $\eta'$ periods, which give charts
near the boundary. We abbreviate $a_{i} = \omega'(\alpha_{i})$ and
$b_{i} = \omega'(\beta_{i})$.
\par  
Next we decompose $\beta_{j} = \beta_{j}^X +  \beta_{j}^\circ$ into the
'eXterior' part $\beta_{j}^X$ outside the plumbing fixture and the
part $\beta_{j}^\circ$ between the two seams of the plumbing fixture,
as in Figure~\ref{cap:betadecomp}.
\input{welded_II}

The separation happens at fixed
sections (of the universal family over the stratum $\omoduli[g-kh,\wh{n}+2kh]
(\wh\mu, (-1)^{2kh})$) in the neighborhood of $(\wX',\eta')$ in the plumbing
annuli~$\cBB_j$, say at the points $u = \delta_0$ and $v = \delta_0$.
Equation~\eqref{eq:horplumbform} simplifies in the one-level case to
$\Omega_{j} = r_{j} dv/v$ where $r_{j} = \zeta^j a_{1}/2\pi i=a_{j}/2\pi i$.
We compute
$$ b_j \= \int_{\beta_{j}} \omega'
\= \int_{\beta_{j}^X} \eta' + \int_{\delta_0}^{x/\delta_0}\Omega_j 
\= \int_{\beta_{j}^X} \eta' + r_{j} (\log x - 2\log \delta_0)\,,
$$
which is well-defined in $\CC +  r_{j}\ZZ$ because of the ambiguity of~$\beta_j$.
By definition of the area form and since $\omega'$ and $\eta'$ agree outside
the plumbing fixture,
\ba \label{eq:areainsympl1}
h(X,q)^{1/k}  &\= \frac{i}{2} \sum_{j=1}^{k} (a_{j} \ol{b_{j}}
- b_{j} \ol{a_{j}})
\,+\, \frac{i}{2} \sum_{j=k+1}^{\wh{g}} (a_j \ol{b_j} - b_j \ol{a_j}) \\
&\= C 
\,+\, \frac{i}{2} \sum_{j=k+1}^{\wh{g}} (a_j \ol{b_j} - b_j \ol{a_j})
-   \frac{k}{4\pi} \cdot |a_1|^2\log(|x|^2)
\ea
is independent of the ambiguity of~$b_{j}$. Here~$C$ is some
function that stems from the integration in the thick part and that is
independent of $x$. We may now let $h^{\text{tck}}_{(0)} = C 
\,+\, \frac{i}{2} \sum_{j=k+1}^{\wh{g}} (a_j \ol{b_j} - b_j \ol{a_j})$
and $R^{\text{hor}}_{(0),1} =  \frac{k}{4\pi} \cdot \pi|a_1|^2$. Both functions are
smooth and bounded away from zero near $(\wX,\bfomega,\wh\Gamma)$,
in fact they correspond to the volume of the region outside the handles and
the residue at the handle respectively. 
\par
For a {\bf general~$X$ that has only horizontal nodes} we
arrive at a similar	formula. We decompose the plumbed surface of the
canonical covering into the thick part and the plumbing fixtures $\WW_j$,
for $j=1,\dots, n(0)$. Since the flat area is additive, we can write it as
a sum of the contribution of the flat area of the thick part and the flat
area of the~$\WW_j$, as we did in the previous case. The area of the thick
part is clearly a smooth function of the period coordinates and bounded
away from zero. For each node~$j$ of~$X$, we get $k$-plumbing
fixtures~$\WW_{j,l}$. Since the residues $r_{j,l}$ of the associated simple
pole differential are $\tau$-conjugates for each fixed $j$, they all have
the same modulus that we denote by $|r_j|$. From the computation in the
plumbing fixtures as in the previous case
we see that the flat area is given by
\begin{equation*}\label{eq:areainsympl2}
h(X,q)^{1/k} \=h^{\text{tck}}_{(0)}-\frac{k}{4\pi}\sum_{j=1}^{E_{(0)}^h} |r_j|^2\log(|x_j|^2)\,
\end{equation*}
which is of the shape we claimed.
\end{proof}
\par 
Now suppose we work in a neighborhood~$U$ of a general boundary point with
notations $\left(\bfx\,;\bft\,;\,\coprod_{i = 0}^{L} \PPer_{(-i)}\right)$ for the perturbed
period coordinates as in Proposition~\ref{prop:PPerIntro} and in detail
in~\eqref{eq:PPerexpl}. More precisely, we group the vector~$\bfx$
of coordinates for opening the horizontal nodes as $\bfx = (x_{(-i),j})$,
where~$-i$ denotes the level that contains the nodes and where
$j=1,\ldots,E_{(i)}^h$ labels these nodes. 
\begin{Prop}\label{prop:metricboundary}
On the neighborhood~$U$ the metric~$h$ has the form 
\begin{equation}\label{eq:hshapegeneral}
h(X,q)^{1/k} \= \sum_{i=0}^L  |\tui|^2 \left(h^{\text{tck}}_{(-i)}
+h^{\text{ver}}_{(-i)}+h^{\text{hor}}_{(-i)}\right)
\end{equation}
where $h^{\text{tck}}_{(-i)}$ are smooth positive functions bounded away from zero and 
\begin{align}\label{eq:cfunctions}
h^{\text{ver}}_{(-i)}:=-\sum_{p=1}^{i}R^{\text{ver}}_{(-i),p} \log\left|t_p\right|,\quad h^{\text{hor}}_{(-i)}:=- \sum_{j=1}^{E_{(-i)}^h} R^{\text{hor}}_{(-i),j} \log|x_{(-i),j}|\,,
\end{align}
with $R^{\text{ver}}_{(-i),p}$  is a smooth non-negative function and $R^{\text{hor}}_{(-i),j}$   is a smooth positive function
bounded  away from zero, both involving only coordinates in $\PPer_{(-i)}$.
\end{Prop}
\par
\begin{proof}
Recall that
the plumbing construction decomposes	the surface along the vertical edges
into various levels (and the	plumbing cylinders between the levels). 
For each level,	we decompose the plumbed surface as the union of the
horizontal plumbing fixtures, the vertical plumbing fixtures and the thick
part. 
\par
 In  \autoref{le:smoothvert}, we have investigated the contribution of  the vertical fixtures, while in \autoref{le:smoothhor} we have investigated the contribution of  the horizontal fixtures. By defining $h^{\text{tck}}_{(-i)}$ the area of the thick part and since the contribution of the horizontal plumbing fixtures at level $-i$ have to be rescaled by $|\tui|^2$, by summing together the contribution of \eqref{eq:metricvertical} and \eqref{eq:metrichor}, we have shown the claim.
\end{proof}
\par
Before proceeding to the proof of Theorem~\ref{thm:areagood}
we recall as an aside and for comparison the definition of a good
metric in the sense of~\cite{mumford77} on a smooth $r$-dimensional
variety (or orbifold)~$\ol{X}$.
\par
Suppose that~$\ol{X}$ is the compactification of~$X$ with a normal
crossing boundary divisor $\partial X = \ol{X} \setminus X$. 
Let~$\cLL$ be a line bundle on~$\ol{X}$. A  metric~$h$ on~$\cLL|_X$ 
is {\em good}, if for each point $p \in \partial X$
there is a neighborhood $\Delta^r$ with coordinates such that
$\partial X = \{\prod_{i=1}^k x_i = 0\}$ and such that the function
$h_s = h(s,s)$ for a local generating section~$s$ of~$\cLL$
has the following properties:
\begin{itemize}
\item[(i)] There exist $C>0$ and $n \in \NN$ such that
$|h_s| < C \Bigl(\sum_{i=1}^k \log|x_i|\Bigr) ^{2n}$ and
$|h_s^{-1}| < C \Bigl(\sum_{i=1}^k \log|x_i|\Bigr) ^{2n}$.
\item[(ii)] the connection one-form $\partial \log h$ and
the curvature two-form $\ol{\partial}\partial \log h$ have
Poincar\'e growth.
\end{itemize}
Here a $p$-form~$\eta$ is said to have {\em Poincar\'e growth} on $\Delta^r$
if for any choice of sections~$v_i$ of $T_{\ol{X}}(\Delta^r)$ there is~$C$
such that 
$$|\eta(v_1,\ldots, v_p)|^2 \,\leq \, C \prod_{i=1}^p\omega_P(v_i,v_i)$$
holds for $\omega_P$ the product of the Poincar\'e metrics
$|dx_i|^2/x_i^2 \log|x_i|^2$ in the coordinates~$x_i$ for $i \leq k$
and the euclidean metric in the other coordinates.
\par
Mumford shows (\cite[Theorem~1.4]{mumford77}) that for a good metric~$h$
the curvature form $\tfrac{i}{2\pi}[F_h]$ defines a closed $(1,1)$-current
that represents the first Chern class of~$\cLL$. This estimate boils
down to the observation that the 'Poincar\'e' integral
\be \label{eq:poincareint}
\frac{1}{2\pi i}\int_{\Delta_\ve} \frac{dx d\bar{x}}{|x|^2 (\log |x|^2)^2} 
 \= - \int_0^\ve \frac{ds}{s (\log(s))^2} \= \frac{1}{\log \ve} < \infty. 
\ee
 of the Poincar\'e metric on the (punctured) disc~$\Delta_\ve$ is finite and
goes to zero as $\ve \to 0$.
\par
The metric~$h$ is indeed good if there is only one level, i.e.\ if the graph
has no vertical edges, as one can deduce from the estimates in the
propositions below. However {\em the metric~$h$ fails to be good}
if there are several levels and horizontal nodes on lower level.
Consider the simplest such case of a graph with
two levels, one vertex at each level and two edges, one edge joining the
levels and a horizontal edge on lower level. Simplifying the situation
by assuming $k=1$ and that $R^{\text{ver}}_{(-1),1}=0$, that $R^{\text{hor}}_{(-1),1}=2$ and that the other bounded functions are $1$, 
the metric is then given by
\bes
h(X,q) \= 1 + |t_1|^2(1  -2\log|x|)\,.
\ees
We observe that this metric is not good in the sense of Mumford 
near the point $(t_1,x)=(0,0)$, when considering the natural
boundary $\{x = 0\} \cup
\{t_1 = 0\}$  consisting of the complement of
the locus where the metric smoothly extends. 
\par
Suppose the metric were good. Then we would have
a constant~$C$ such that
\bes
\Bigl|\partial \log h\Bigl(\frac{\partial}
{\partial t_1}\Bigr)\Bigr|^2 \leq \frac{C}{|t_1|^2 (\log|t_1|)^2}
\ees
on the neighborhood~$U$ of $(t_1,x)=(0,0)$, which is
equivalent to the inequality of the square roots
\be \label{eq:sqrtgood}
\frac{|t_1| (1-2\log|x|)}{1 + |t_1|^2(1-2\log|x|)}
\leq \frac{C^{1/2}}{|t_1| \log|t_1|}
\ee
Choosing a sequence tending to $(0,0)$ with $1-\log|x| = |t_1|^{-2}$
we get a contradiction.
\par
Instead of aiming for a bound as in the definition of good, integrability
statements are sufficient. In fact the
coefficients of
\bes
\partial \log(h) = \frac{|t_1|^2 (1-2\log|x|)}{h}  \frac{dt_1}{t_1}
- \frac{|t_1|^2}{h} \frac{dx}x
\ees
and
\begin{align*}
\ol{\partial} \partial \log h
&= \frac{|t_1|^2(1-2\log|x|)}{h^2} \frac{d\bar{t}_1}{\bar{t}_1}\frac{dt_1}{t_1}
-\frac{|t_1|^4}{h^2} \frac{d\bar{x}}{\bar{x}} \frac{dx}x-\frac{|t_1|^2}{h^2} \left(\frac{d\bar{t}_1}{\bar{t}_1}\frac{d{x}}{x}+ \frac{d\bar{x}}{\bar{x}}\frac{dt_1}{t_1}\right)
\end{align*}
are locally integrable, and thus define currents. To see that the current
$F_h = [\ol{\partial} \partial \log h]$ given by the curvature form is
closed, we have to show that we can apply the derivative (in the sense
of currents) inside the brackets, on the differential form, where it gives
zero. This requires an application of Stokes' theorem, and thus an
integral over the boundary $T_\delta$ of a shrinking tubular neighborhood
around the locus where
the metric is not smooth, i.e.\ around $\{x=0\}$ and $\{t_1=0\}$. To see that this current
represents the first chern class $c_1(\cOO(-1))$, we compare with the curvature
form of a smooth metric. To see that the difference is zero in cohomology,
the term $[d\partial \log(h)]$ appears and we'd like to invoke say that
this is $d[\partial \log(h)]$, i.e.\ a coboundary of a current. This is
gives a second application of Stokes' theorem, justified by another
integration over~$T_\delta$. To justify that we can pass to wedge powers
is a third application of  Stokes' theorem. This integrals are estimated
in the general case in the following proofs.
\par
\medskip
The {\em Proof of Theorem~\ref{thm:areagood}} is now contained in the
following two propositions.
\par
\begin{Prop} \label{prop:arecurrents}
The differential forms $\Omega = \partial \log h$ and
$F_h = \ol{\partial}\partial \log h$, and more generally the forms
$F_h^d$ and $\Omega \wedge F_h^d$ have coefficients in $L^1_{\rm loc}$.
\par
In particular $F_h^d$ defines a current of type $(d,d)$ for any~$d \in \NN$.
\end{Prop}
\par
\begin{Prop} \label{prop:represents}
The current $[\ol{\partial} \partial \log h]$ is closed and $\tfrac{1}{2\pi i}$
times the curvature $(1,1)$-form $F_h = \ol{\partial} \partial \log h$
represents the first Chern class $c_1(\cOO(-1))$ in cohomology.
\par
More generally, the wedge powers $\wedge^d \Bigl(\tfrac{1}{2\pi i} F_h\Bigr)$
represent the class of $c_1(\cOO(-1))^d$ in cohomology.
\end{Prop}
\par
We calculate the relevant differential forms explicitly. More specifically, we
first determine explicitly the types of differential forms
that we can encounter in $\d h$, $\ol{\d} h$ and $\ol{\d}\d h$ up to
continuous factors that don't affect integrability.
\par
Recall that by \autoref{prop:metricboundary}, the function $h$ near a boundary point is given as the sum of three contributions given by the thick part and the vertical and horizontal plumbing fixtures. 
\par
For every level $(-i)$, the thick part contribution $ \sum_{i=0}^L |\tui|^2 h^{\text{tck}}_{(-i)}$ is a smooth positive function bounded away from zero. Hence all of its derivatives are in particular smooth.
\par
Then we analyze the contribution $ \sum_{i=0}^L |\tui|^2 h^{\text{ver}}_{(-i)}$  from the vertical plumbing fixtures. An important remark is that the functions $|\tui|^2h^{\text{ver}}_{(-i)}$ are continuous, since they are given by sums of functions  of type $|t_p|^2\log\left|t_p\right|$, for $p\leq i$. By using the explicit expression of $h^{\text{ver}}_{(-i)}$ given in \eqref{eq:cfunctions} we compute
\begin{align*}
	\d\left(  |\tui|^2 h^{\text{ver}}_{(-i)}\right)=&-|\tui|^2\sum_{p,p_1=1}^i  \ell_{p_1}R^{\text{ver}}_{(-i),p} \log\left|t_p\right| \frac{dt_{p_1}}{t_{p_1}}\\
	&\ -|\tui|^2\sum_{p=1}^i \left( \frac{1}{2}R^{\text{ver}}_{(-i),p} \frac{dt_{p}}{t_{p}}+\log\left|t_p\right|\cdot \d R^{\text{ver}}_{(-i),p} \right).
\end{align*}
It is clear that all the  1-forms appearing in the previous expression are continuous. The same is obviously true for the analogous form given by the $\ol{\d}$ derivative.
We compute now the second derivative 
\begin{align*}
	&\ol{\d}\d\left(  |\tui|^2 h^{\text{ver}}_{(-i)}\right)=
	-|\tui|^2\sum_{p,p_1,p_2=1}^i \ell_{p_1} \ell_{p_2} R^{\text{ver}}_{(-i),p} \log\left|t_p\right|
	\frac{d\bar{t}_{p_2}}{\bar{t}_{p_2}}
	\frac{d{t}_{p_1}}{{t}_{p_1}} \\
	&\quad \quad\quad-|\tui|^2\sum_{p,p_1=1}^i \ell_{p_1}\left( \frac{1}{2}R^{\text{ver}}_{(-i),p}\frac{d\bar{t}_{p_1}}{\bar{t}_{p_1}} \frac{dt_{p}}{t_{p}}+\log\left|t_p\right|\frac{d\bar{t}_{p_1}}{\bar{t}_{p_1}}\wedge \d R^{\text{ver}}_{(-i),p} \right)\\
	&\quad \quad\quad -|\tui|^2\sum_{p,p_1=1}^i \ell_{p_1}\left( \log\left|t_p\right| \ol{\d}R^{\text{ver}}_{(-i),p}\wedge \frac{dt_{p_1}}{t_{p_1}}+\frac{1}{2} R^{\text{ver}}_{(-i),p}\frac{d\bar{t}_{p_1}}{\bar{t}_{p_1}}
	\frac{d{t}_{p}}{{t}_{p}}\right)\\
	&\quad \quad\quad -|\tui|^2\sum_{p=1}^i	\left( \frac{1}{2}\ol{\d} R^{\text{ver}}_{(-i),p} \wedge\frac{dt_{p}}{t_{p}}+\frac{1}{2}\frac{d\bar{t}_{p}}{\bar{t}_{p}} \wedge\d  R^{\text{ver}}_{(-i),p}+\log\left|t_p\right|\cdot \ol{\d}\d R^{\text{ver}}_{(-i),p} \right).
\end{align*}
By inspecting the terms of the previous expression, we can see that the only
non-continuous coefficients that can appear stem from the first line. Indeed
if $p=p_1=p_2$ and $\ell_{p}=1$, then we can have 2-forms (up to multiplication
by smooth positive functions) of type 
\begin{equation}\label{eq:vertterms}
	\log\left|t_{p}\right| d\bar{t}_{p}d{t}_{p}.
\end{equation}
\par
We finally analyze the contribution $ \sum_{i=0}^L |\tui|^2 h^{\text{hor}}_{(-i)}$  from
the horizontal plumbing fixtures. Again by using the explicit expression
of $h^{\text{hor}}_{(-i)}$ given in \eqref{eq:cfunctions} we compute
\begin{align*}
	\d\left(  |\tui|^2 h^{\text{hor}}_{(-i)}\right)=&- |\tui|^2 \sum_{p=1}^i\sum_{j=1}^{E_{(-i)}^h} R^{\text{hor}}_{(-i),j} \log|x_{(-i),j}| \frac{dt_p}{t_p}\\
	&\ - |\tui|^2 \sum_{j=1}^{E_{(-i)}^h} \left(\log|x_{(-i),j}|\cdot  \d R^{\text{hor}}_{(-i),j}+\frac{1}{2}R^{\text{hor}}_{(-i),j} \frac{dx_{(-i),j}}{x_{(-i),j}} \right).
\end{align*}
The non-continuous 1-forms appearing in the previous expression (up to continuous coefficients) are
\begin{equation}\label{eq:hor1terms}
	 |\tui|^2 \log|x_{(-i),j}| \frac{dt_p}{t_p},\quad |\tui|^2 \log|x_{(-i),j}| \d R^{\text{hor}}_{(-i),j},\quad |\tui|^2  \frac{dx_{(-i),j}}{x_{(-i),j}}
\end{equation}
for $p\leq i$.
We compute now the second derivative
\begin{align*}
	&\ol{\d}\d\left(  |\tui|^2 h^{\text{hor}}_{(-i)}\right)=- |\tui|^2 \sum_{p_1,p_2=1}^i\sum_{j=1}^{E_{(-i)}^h} R^{\text{hor}}_{(-i),j} \log|x_{(-i),j}| \frac{d\bar{t}_{p_2}}{\bar{t}_{p_2}} \frac{dt_{p_1}}{t_{p_1}}\\
	&\quad\quad\quad\ - |\tui|^2 \sum_{p=1}^i\sum_{j=1}^{E_{(-i)}^h} \left(\log|x_{(-i),j}|\frac{d\bar{t}_{p}}{\bar{t}_{p}} \wedge  \d R^{\text{hor}}_{(-i),j}+\frac{1}{2}R^{\text{hor}}_{(-i),j} \frac{d\bar{t}_{p}}{\bar{t}_{p}} \wedge\frac{dx_{(-i),j}}{x_{(-i),j}} \right)\\
	&\quad\quad\quad- |\tui|^2 \sum_{p=1}^i\sum_{j=1}^{E_{(-i)}^h} \left( \log|x_{(-i),j}| \ol{\d}R^{\text{hor}}_{(-i),j}\wedge \frac{dt_p}{t_p}+\frac{1}{2}R^{\text{hor}}_{(-i),j}\frac{d\bar{x}_{(-i),j}}{\bar{x}_{(-i),j}}  \wedge \frac{dt_p}{t_p}\right)\\
	&\quad\quad\quad\ - |\tui|^2 \sum_{j=1}^{n(-i)}\left( \frac{1}{2}\frac{d\bar{x}_{(-i),j}}{\bar{x}_{(-i),j}}  \wedge  \d R^{\text{hor}}_{(-i),j}+ \log|x_{(-i),j}| \cdot \ol{\d} \d R^{\text{hor}}_{(-i),j}+ \frac{1}{2}\ol{\d}R^{\text{hor}}_{(-i),j} \wedge \frac{dx_{(-i),j}}{x_{(-i),j}} \right).
\end{align*}
The non-continuous 2-forms appearing in the previous expression (up to continuous coefficients) are
\begin{align}\label{eq:hor2terms}
	&|\tui|^2\log|x_{(-i),j}| \frac{d\bar{t}_{p_2}}{\bar{t}_{p_2}} \frac{dt_{p_1}}{t_{p_1}},\quad	|\tui|^2 \log|x_{(-i),j}|\d R^{\text{hor}}_{(-i),j}   \frac{d\bar{t}_{p}}{\bar{t}_{p}} ,\quad  |\tui|^2 \frac{d\bar{t}_{p}}{\bar{t}_{p}} \frac{dx_{(-i),j}}{x_{(-i),j}}  \\
	&\nonumber |\tui|^2 \log|x_{(-i),j}|  \ol{\d} \d R^{\text{hor}}_{(-i),j},\quad |\tui|^2\d R^{\text{hor}}_{(-i),j} \frac{d\bar{x}_{(-i),j}}{\bar{x}_{(-i),j}}
\end{align}
and their complex conjugates, where $p,p_{1},p_{2}\leq i$. Recall also that by \autoref{prop:metricboundary} the functions   $R^{\text{hor}}_{(-i),j}$   are smooth and involve only coordinates in $\PPer_{(-i)}$.
\par 
In order to analyze the general expression for the differential forms
associated with the metric~$h$, note that 
\begin{align}\label{eq:dlog}
\partial \log(h) =  k\frac{\partial h^{1/k}}{h^{1/k}},\quad\quad
F_h \= \ol{\partial} \partial \log h &\=  k\frac{\ol{\partial} \partial h^{1/k}}{h^{1/k}}
\,-\, k\frac{\ol{\partial} h^{1/k}}{h^{1/k}} \wedge \frac{\partial h^{1/k} }{h^{1/k}}.
\end{align}
\par
Hence the non-continuous 1-forms appearing in $\partial \log(h)$ consist of
linear combinations of the building blocks given by the quotient by $h^{1/k}$
of the terms
appearing in~\eqref{eq:hor1terms}, while the non-continuous 2-forms appearing
in $\ol{\partial} \partial \log(h)$ consist of the two-fold wedge products
of type $(1,1)$ of the quotient by $h$ of the terms appearing in \eqref{eq:hor1terms}
and of the terms appearing in~\eqref{eq:vertterms} and~\eqref{eq:hor2terms},
together with their complex conjugates. 
\par
We now fix some more notation. We may assume that our neighborhood~$U$ is the
product of the polydiscs
\be
D^\bft \= \{\bft : t_i\in \Delta_\ve\} \quad \text{and} \quad
D^\bfx \= \{\bfx : x_{(-i),j} \in \Delta_\ve\}\,.
\ee
in the corresponding variables times a ball~$B$ corresponding to all the
variables in $\PPer_{(-i)}$ for $i=0,\ldots,L$.
\par
\begin{proof}[Proof of Proposition~\ref{prop:arecurrents}]
We first record the following (``sharp'') bound for the reciprocal of the metric
\begin{equation}\label{eq:areabound2}
\frac{1}{h^{1/k}}\leq  \frac{1}{h^{\text{tck}}_{(0)} +|\tui|^2 h^{\text{hor}}_{(-i)}}\leq
\frac{1}{C-|\tui|^2\log|x_{(-i),j}|}, 
\end{equation}
for some constant $C$, that stems from the fact that
$h^{\text{tck}}_{(0)}$ is uniformly bounded away from zero on $U$.
This obviously implies that $\frac{1}{h^{1/k}}$ is bounded, which we refer
to as the "first weak bound" and that$ \frac{1}{h^{1/k}}
\leq -1/|\tui|^2\log|x_{(-i),j}|$, which we call the "second weak bound".
Hence, using the expressions~\eqref{eq:dlog}, all terms appearing in $\d h^{1/k}$ and
$\ol{\d} \d h^{1/k}$ with continuous coefficients induce locally integrable
forms in $\partial \log h$, $\ol{\partial} \partial \log h$ and their powers.
\par
We now treat the non-continuous terms appearing in $\partial \log h$ and
$\ol{\partial} \partial \log h$,  i.e. the quotient by $h^{1/k}$ of the terms
displayed in  \eqref{eq:hor1terms} and \eqref{eq:vertterms},  \eqref{eq:hor2terms}
respectively. We refer to these terms as 'building blocks'.
\par
We first look at each term in $\partial \log h$, i.e. we focus on the quotient
by~$h^{1/k}$ of the coefficients~\eqref{eq:hor1terms}.  For $dt_{p}$-differential forms
we use the weak bound for $1/h^{1/k}$  and the fact that  $|\tui|^2/t_p$ is a
polynomial expression in the $t_i$ and $\bar{t}_i$ for $p\geq i$. The logarithmic
contribution is unbounded, but after a change to polar coordinates we are left with
\ba \label{eq:dxxint}
& \phantom{\=} \int_{B \times D^\bfx \times D^\bft} \frac{|\tui|^2}{h^{1/k} \cdot t_p}
\log|x_{(-i),j}|  d{\rm vol}
\leq C_1 \int_{ \Delta_\ve}\log|x_{(-i),j}||dx_{(-i),j}|^2 \\
& \leq C_2\int_{r=0}^{\ve}
r \log |r|dr < \infty\,.
\ea
The same argument applies to the $|\tui|^2 \log|x_{(-i),j}| \d R^{\text{hor}}_{(-i),j}$ term.
For the $dx_{(-i),j}$-coefficient we use again the weak bound for $1/h^{1/k}$ and so we
are left with a polynomial expression in the $t_i$ and $\bar{t}_i$ and
the finite integral 
\begin{equation}\label{eq:dxoverxint}
\int_{\Delta_\ve} \frac{1}{x_{(-i),j}}|dx_{(-i),j}|^2 <\infty.
\end{equation}
The same arguments apply verbatim to the coefficients of $\ol{\partial} \log h$.
\par
Next, we examine the coefficients in $\frac1{h^{1/k}} \ol{\partial} \partial h^{1/k}$.
Using the weak bound  and the same polar coordinates argument as
in~\eqref{eq:dxxint}, we see that the term $\log\left|t_{p}\right| |dt_{p}|^2/h^{1/k}$
appearing in \eqref{eq:vertterms} is locally integrable. The other terms given by
the quotient by $h^{1/k}$ of the terms appearing in \eqref{eq:hor2terms} are of the same shape as the one already treated, so they are locally integrable.
\par
Finally we examine the terms that may arise from as an arbitrary wedge
product of $\tfrac1{h^{1/k}} \ol{\partial} \partial h^{1/k}$,
$\tfrac1{h^{1/k}} \partial h^{1/k}$ or $\tfrac1{h^{1/k}}\ol{\partial} h^{1/k}$.
This relies on Fubini and a suitable organization of the order of integration,
necessary since the building blocks with $dt_p$ may have coefficients involving
$x_{(-i),j}$ and vice versa. Ultimately we rely on the observation that in a non-zero
wedge product each of the $dt_p$ and $dx_{(-i),j}$ and their conjugates appear only once.
\par
To start, observe that the building block~\eqref{eq:vertterms} has no
$x_{(-i),j}$-dependence. Whenever such a term appear at a level below any level
that has a horizontal node, we integrate $|dt_p|^2$ and use the integrability
of $\log|t_p|$. In general, for any level $(-i)$ and a subset $P_i\subseteq \{1,\ldots,i\}$, we will use the notation 
\[f_{P_i}(\bft ):=\prod_{p \in P_i}|t_p|\log|t_p|\]
to indicate the functions appearing in the product of the building blocks of type\eqref{eq:vertterms}  $	\prod_{p  \in P_i}\log|t_p||dt_p|^2 = f_{P_i}(\bft ) \prod_{p \in P_i} d\theta_pd|t_p|$
after passing to polar coordinates polar coordinates. We record for the sequel that the functions $f_i(\bft )$ are continuous in a neighborhood of zero.
\par
Whenever we have a combination of the forms in~\eqref{eq:hor1terms} and
in~\eqref{eq:hor2terms}  not containing the $dx_{(-i),j}$, we can use the first
weak bound for $1/h^{1/k}$. Then we are left with  $\log|x_{(-i),j}|^N$ for
some $N>0$, which is integrable.
\par
If a differential form stems from a derivative of
$R^{\text{hor}}_{(-i),j}$, the ratio of the coefficient over $h^{1/k}$  can be bounded
by a constant using~\eqref{eq:areabound2} (or rather just the second weak bound).
We can thus disregard those  $R^{\text{hor}}_{(-i),j}$-derivatives in the rest of
the discussion.
\par
We now treat the case where we have a combination of building blocks involving
$dx_{(-i),j}$, i.e., coming from the first term of~\eqref{eq:hor1terms} and its
conjugate, possibly together with other building blocks whose coefficients
involve $\log(x_{(-i),j})$. 
\par
\medskip
\noindent
\textbf{Case A: Only one of $dx_{(-i),j}$ and its conjugate appear from
building blocks.} We perform the integral over $|dx_{(-i),j}|^2$ first.
This leads to the integral as in~\eqref{eq:dxoverxint}, or maybe with an
additional factor $\log(x_{(-i),j})^N$, that does not change the finiteness
of the integral. Using the first weak bound for $h^{1/k}$ a polynomial expression
in the $t_i, \ol{t}_i$ remains, which is irrelevant to the finiteness discussion
for the subsequent integrals.
\par
\medskip
\noindent
\textbf{Case B: Both $dx_{(-i),j}$ and its conjugate appear, but with
no other building block involving $\log(x_{(-i),j})$}
Then the second weak estimate yields the same situation as the
Poincar\'e integral~\eqref{eq:poincareint}, which is finite. We perform
these integrals before addressing the other $dx_{(-i),j}$. For simplicity
of notation we label the $x$-variables subsequently by $x_e$ for $e \in E^h$,
keeping in mind that such a horizontal edges comes with a level $i=i(e)$.
We denote the set of horizontal edges that do not belong to Case~$A$ or~$B$
by $E^h_C \subset E^h$.
\par
\medskip
\noindent
\textbf{Case C: Both $dx_{(-i),j}$ and its conjugate appear and there 
are building blocks  involving $\log(x_{(-i),j})$} We denote by
$i_\max := \max_{e\in E_C^h} |i(e)|$ the largest index of a level that has
a horizontal node. 
Define $P_e^1 \subset I:= \{1,\ldots,i_\max\}$ to be the set of indices where
the product of building blocks $\frac{|\tui|^4\log|x_{e}|^2}{h^{2/k}}
\frac{|dt_{p}|^2}{|t_{p}|^2}$ (coming from the last term of~\eqref{eq:hor1terms}
and its conjugate) occur.
Define $P_e^2 \subset I$ the set of indices parametrizing forms of
type $\frac{|\tui|^2\log|x_{e}|}{h^{1/k}} \frac{|dt_{p}|^2}{|t_{p}|^2}$ (coming
from the first term of~\eqref{eq:hor2terms}). We define $P_e^3 \subset I$, resp.\
$P_e^4 \subset I$ the indices, where last term of~\eqref{eq:hor1terms} but not
its conjugate (resp. the other way round) occur.
Note that in any non-zero building block the subsets $P_e^j$ are disjoint
and that the subsets $P_e = P_e^1 \cup P_e^2 \cup P_e^3 \cup P_e^4$ 
are disjoint for different~$e$, since otherwise the wedge products
of building blocks is zero. 
Wedging all these building blocks together we obtain the
differential form 
\be\label{eq:badcaseP}
\prod_{e \in E^h_C} 
\left(\frac{|t_{\lceil i(e)\rceil}|^{2(2 + M_e)}} {h^{\frac{2+M_e}{k}}} 	(\log|x_{e}|)^{M_e}
\prod_{p \in P_e^1  \cup P_e^2 } \frac{|dt_p|^2}{|t_p|^{2}} \prod_{p \in P_e^3}
\frac{|dt_p|^2}{t_p}\prod_{p \in P_e^4} \frac{|d\ol{t}_p|^2}{\ol{t}_p} \right)
\frac{|dx_{e}|^2}{|x_{e}|^2}\,,
\ee
where $M_e = 2|P_e^1| +|P_e^2\cup P_e^3 \cup P_e^4|$. We use that 
$\prod_{p \in P_e} t_p$ divides $t_{\lceil i(e)\rceil}$ to cancel the $t_p$ and $\ol{t}_p$
in the denominator and apply $M_e$ times the
second weak estimate in each factor to cancel all terms raised to~$M_e$.
Note that the expression~\eqref{eq:badcaseP} involves all variables $t_i$
for $i \in I$. We thus integrate also over $\prod_{p \in I \setminus P_e}$
and have to take into account the differential form that stems from
the building blocks in~\eqref{eq:vertterms} in the total wedge, which will contribute with the continuous function $f_{P_{i_\max}}(\bft )$ for some $P_{i_\max}\subseteq I$. 
We pass to polar coordinates and write $R_e = |x_{e}|$ and $r_p = |t_p|$.  
Coarsely estimating the angle integral it suffices to show the finiteness
of the following expression:
\begin{align}\label{eq:badintegral2}
&\int_{0}^{\ve} \cdots \int_{0}^{\ve} f_{P_{i_\max}}(\bft ) \cdot
\prod_{e \in E^h_C}  \frac{\prod_{p\leq i(e)} r_{p}^3}{\left(C- \log(R_e)
\prod_{p\leq i(e)} r_{p}^2\right)^2} \frac{dR_e}{R_e} \prod_{i\in I} dr_{i}  \\
&\nonumber =\int_{0}^{\ve } \cdots  \int_{0}^{\ve}f_{P_{i_\max}}(\bft )\cdot
\left( \prod_{e \in E^h_C}  \left(\prod_{p\leq i(e)} r_{p}
\int_{-\infty}^{\log(\ve)\prod_{p\geq i(e)} r_{p}^2} \frac{  dv_e }{(C-v_e)^2}
\right)\prod_{i\in I} dr_{i} \right) \\
&\nonumber =\int_{0}^{\ve}  \cdots  \int_{0}^{\ve} f_{P_{i_\max}}(\bft ) \cdot \prod_{e \in E^h_C} 
\left(	\frac{ \prod_{p\leq i(e)} r_{p} }{C-\log(\ve)\prod_{p\leq i(e)} r_{p} ^2}
\right)\prod_{i\in I} dr_{i}<\infty
\end{align}
where we have used the change of coordinates $v_e=\log(R_e)\prod_{p\geq i(e)} r_{p}^2$
in the first equality. The last expression is finite since we integrate the product of continuous functions.
\par
Since we have shown that all the forms $F_h^d$ and $\Omega \wedge F_h^d$ have
locally integrable coefficients, their corresponding currents of integration
are well-defined.
\end{proof}
\par
\begin{proof}[Proof of Proposition~\ref{prop:represents}]
We identify the local statements needed to prove the
claims and justify them subsequently.
To see that $[F_h] = [\ol{\partial} \partial \log h]$ defines a closed
current we need to justify the first step in the chain
\be
d[\ol{\partial} \partial \log h] \= [d(\ol{\partial} \partial \log h)] \= 0
\ee
of cohomology classes of currents. By definition we have to justify that
\be
\int_{D \times B} dF_h \wedge \xi = - \int_{D \times B} F_h \wedge d\xi
\ee
for any smooth $r$-form~$\zeta$, where $r = \dim_\RR \PP \kLMS-3$. By
Stokes' theorem this amounts to justify that  
\be \label{eq:tub1}
\lim_{\delta \to 0} \int_{T_\delta} F_h \wedge \xi \= 0
\ee
where $T_\delta$ is one of the tubular neighborhoods inside $B \times D^\bft \times D^\bfx$, with tube radius~$\delta$, around the divisors defined by setting one coordinate axis to zero. We will denote such tubular neighborhoods by 
\bes
T^{(-i),j}_\delta \= B \times \Bigl\{|x_{(-i),j}| = \delta \,;\,\, t_{p}
\in \Delta_\ve\,\,
\text{for all $p$}\,;\,\,  x_{(-i'),j'} \in \Delta_\ve \,\,
\text{for all $(i',j') \neq (i,j)$}
\Bigr\}
\ees
and
\bes
T^{p}_\delta \= B \times \Bigl\{ t_{p}= \delta \,;\,\, t_{p'}
\in \Delta_\ve\,\,
\text{for all $p'\neq p$}\,;\,\,  x_{(-i),j} \in \Delta_\ve \,\,
\text{for all $ (i,j)$}
\Bigr\}
\ees
and finally by
\bes
T_\delta^B=D\times \partial B^\delta,
\ees
where $\partial B^\delta $ is the union of the tubular neighborhoods
around the coordinate axis in the $\PPer_{(-i)}$ components.
\par
For the second statement of the statement let $h^*$ be a smooth (comparison)
metric on  $\cOO(-1)$. Then certainly $\tfrac{1}{2\pi i}$ times the curvature
$F_h^* = \ol{\partial} \partial \log h^*$ represents the first
Chern class of $\cOO(-1)$. To justify the equality of cohomology classes of
currents 
\bes
[\ol{\partial} \partial \log h^*] - [\ol{\partial} \partial \log h]
\= [d(\partial \log h^* -\partial \log h)]  \= d[\partial \log h^* -\partial \log h]
\= 0
\ees 
we have to justify the second equality sign, i.e. that the current of i
ntegration of $\ol{\partial} \partial \log h$ is the same as the derivative in
the sense of currents of $\partial \log h$. Then the last equality follows from
Proposition~\ref{prop:arecurrents}, showing that the expression is
a coboundary in the sense of currents, since $\log h^* - \log h = \log(h^*/h)$
is independent of the scale of~$h$ and thus globally well-defined.
\par
Writing~$\Omega^* = \partial \log h^*$ and $\Omega = \partial \log h$, we have
to justify that for any smooth $\dim_\RR \PP \kLMS-2$-form the condition
\[\lim_{\delta \to 0} \int_{T_\delta} (\Omega^* - \Omega) \wedge \xi \= 0\] 
holds, where $T_\delta\in\{T^{(-i),j}_\delta,\ T_\delta^{p},\ T_\delta^B\}$ for all $(i,j)$
and~$p$. This follows once we have shown 
\be \label{eq:tub2}
\lim_{\delta \to 0} \int_{T_\delta} \Omega \wedge \xi \= 0
\ee
and from the smoothness of~$\Omega^*$.
\par
For the generalization to wedge powers we use $F_h = d\Omega$ and $F_h^*
= d \Omega^*$ and want to argue that there is
an equality of cohomology classes of currents
\bes
[F_h^d] - [(F_h^*)^d] \= d\Bigl[
(\Omega - \Omega^*) \wedge \sum_{i+j=n-1}F_h^{i}(F_h^*)^{j}
\Bigr]\,.
\ees
With this equation at hand we use that the argument of the differential
operator on the right hand side defines a current by
Proposition~\ref{prop:arecurrents}, so that $[F_h^d]$ and $[(F_h^*)^d]$
are cohomologous and $(\tfrac{1}{2\pi})^d[(F_h^*)^d]$ is known
to represent $c_1(\cOO(-1))^d$.
\par
To justify this equation we need again to argue about the interchange of derivative and passage to the current. Hence we  need to show that  for all~$d$
and for all smooth $\dim_\RR \PP \kLMS-2d-1$-forms $\xi$, the following equation
\be \label{eq:tub3}
\lim_{\delta \to 0} \int_{T_\delta} \Omega \wedge F_h^n \wedge \xi \= 0
\ee
holds for $T_\delta\in\{T^{(-i),j}_\delta,\ T_\delta^{p},\ T_\delta^B\}$ for all $(i,j)$ and $p$.
\par
To justify these three equations~\eqref{eq:tub1}, \eqref{eq:tub2}
and~\eqref{eq:tub3} we fix the tubular neighborhood $T_\delta$ around
one of the boundary divisors and analyze the forms that may appear from the
wedge products in~\eqref{eq:tub1}, in~\eqref{eq:tub2} or in~\eqref{eq:tub3}.
These are wedge products of the quotient by $h$ of the building blocks in~\eqref{eq:vertterms}, \eqref{eq:hor1terms}, \eqref{eq:hor2terms} and the differentials of the coordinates themselves.
\par
Our strategy is to apply Fubini and use \autoref{prop:arecurrents} in order to ensure that the blocks not depending on the coordinate defining $T_\delta$ give a finite result, which is independent of $\delta$. Then we need to check that the integral of arbitrary products of the building blocks appearing in \eqref{eq:tub1}, ~\eqref{eq:tub2} and in~\eqref{eq:tub3} involving the special coordinate, yields an expression going to zero for $\delta$ going to zero.
\par
First of all, we consider the tubular neighborhoods $T_\delta^B$ around the zero divisor  of a coordinate function belonging to the $\PPer_{(-i)}$ part. Since the differential forms appearing as integrands in the three equations~\eqref{eq:tub1}, \eqref{eq:tub2}
and~\eqref{eq:tub3} are continuous in the $\PPer_{(-i)}$ coordinates, by \autoref{prop:arecurrents} we conclude that the integrand over $T_\delta^B$ can be bounded by a constant function, so the integral goes to zero for $\delta$ going to zero.
\par
Next, note that whenever we consider the integral over $T^{(-i),j}_\delta$, resp.\
over $T_\delta^{p}$, of a product of building blocks not depending on the
variables $x_{(-i),j}$, resp.\ on the  $t_{p}$,  then \autoref{prop:arecurrents} yields
finiteness independently of~$\delta$.  
\par
Hence we are left to consider integrals of product of building blocks corresponding
to the variables~$x_{(-i),j}$ and~$t_{p}$. We observe that the differential forms
involving $|dx_{(-i),j}|^2$ resp.\ $|dt_p|^2$ restrict to zero on $T_\delta^{(-i),j}$
resp.\ on $T_\delta^{p}$.  
\par
We consider first  the integral over $T_\delta^{p}$. By the previous remark,  the block \eqref{eq:vertterms} given by $\log|t_p| |dt_p|^2$ restricts to zero on $T_\delta^{p}$. By inspecting the proof of  \autoref{prop:arecurrents}, we see that the only possibly problematic case is the one corresponding to Case C. We have to consider an expression analogous to \eqref{eq:badcaseP}, but in this situation we can only have a $dt_p$ or a $d\bar{t}_p$. Using the same strategy as for \eqref{eq:badcaseP}, we consider the expression analogous to \eqref{eq:badintegral2} obtained after passing to polar coordinates and estimating the angle integral. In this situation we hence have that the integral over $T_\delta^{p}$ of these product of building blocks is given by 
\begin{align*}
	&\int_{0}^{\ve} \cdots \int_{0}^{\ve} f_{P_{i_\max}}(\bft ) \cdot
	\prod_{e \in E^h_C}  \frac{\delta^3\prod_{i\leq i(e), i\not = p} r_{i}^3}{\left(C- \delta^2\log(R_e)
		\prod_{i\leq i(e), i\not = p} r_{i}^2\right)^2} \frac{dR_e}{R_e} \prod_{i\in I\setminus\{p\}} dr_{i}  \\
	&\nonumber \= \int_{0}^{\ve}  \cdots  \int_{0}^{\ve} f_{P_{i_\max}}(\bft ) \cdot \prod_{e \in E^h_C} 
	\left(	\frac{ \delta\prod_{i\leq i(e), i\not = p} r_{i} }{C-\log(\ve)\delta^2\prod_{i\leq i(e), i\not = p} r_{i}^2}
	\right)\prod_{i\in I\setminus\{p\}} dr_{i}\,.
\end{align*}
By estimating by a constant the denominator appearing in the second line, we see that the expression is going to zero for $\delta$ going to zero.
\par
We consider finally the integral over $T^{(-i),j}_\delta$. We already remarked that the
form $|dx_{(-i),j}|^2$ restricts to zero to $T^{(-i),j}_\delta$. We now describe the general case
of a product of building blocks involving the $x_{(-i),j}$ variable. We consider the case
where $dx_{(-i),j}$ appears, the other case where $d\bar{x}_{(-i),j}$ appears is clearly
equivalent. Such a product of building block is given the product of the
expression~\eqref{eq:badcaseP}  (for edges $e$ not corresponding to the special
index $(-i),j$) with
\begin{align*}
\frac{|t_{\lceil i\rceil}|^{2(1 + M_{(-i),j})}} {h^{\frac{1+M_{(-i),j}}{k}}} 	(\log|x_{(-i),j}|)^{M_{(-i),j}}
\prod_{p \in P^1  \cup P^2 } \frac{|dt_p|^2}{|t_p|^{2}} \prod_{p \in P^3}
\frac{|dt_p|^2}{t_p}\prod_{p \in P^4} \frac{|d\ol{t}_p|^2}{\ol{t}_p} 
\frac{dx_{(-i),j}}{x_{(-i),j}}.
\end{align*}
where we dropped the index $(-i),j$ of the sets $P^\ell$.
By proceeding as we did to reach the expression~\eqref{eq:badintegral2} after estimating and
passing to polar coordinates, we obtain that the integral over $T^{(-i),j}_\delta$ of these
product of building blocks is given by the following expression
\begin{align*}
	&\int_{0}^{\ve} \cdots \int_{0}^{\ve} f_{P_{i_\max}}(\bft ) \cdot
	\frac{\prod_{p\leq i} r_{p}}{C- \log(\delta)
		\prod_{p\leq i} r_{p}^2} \prod_{e \in E^h_C}  \frac{\prod_{p\leq i(e)} r_{p}^3}{\left(C- \log(R_e)
		\prod_{p\leq i(e)} r_{p}^2\right)^2} \frac{dR_e}{R_e} \prod_{p\in I} dr_{p}  \\
	&\nonumber =\int_{0}^{\ve}  \cdots  \int_{0}^{\ve} f_{P_{i_\max}}(\bft ) \cdot\frac{\prod_{p\leq i} r_{p}}{C- \log(\delta)
		\prod_{p\leq i} r_{p}^2}  \cdot  \prod_{e \in E^h_C} 
	\left(	\frac{ \prod_{p\leq i(e)} r_{p} }{C-\log(\ve)\prod_{p\leq i(e)} r_{p} ^2}
	\right)\prod_{p\in I} dr_{p}\\
	&\leq C_1\int_{0}^{\ve}  \cdots  \int_{0}^{\ve} \frac{\prod_{p\leq i} r_{p}}{C- \log(\delta)
		\prod_{p\leq i} r_{p}^2} \prod_{p=1}^i dr_{p}\,. \\
\end{align*}
The first inequality is obtained by the same substitutions as in~\eqref{eq:badintegral2}
and the last inequality is given by estimating the product of the continuous functions
by a constant $C_1$.
\par
We will use notation $r_{\lceil i \rceil}:=\prod_{p\leq i} r_{p}$.
If $i=0$, then the last line above is simply given by $C_1/(C-\log(\delta))$ which clearly tends to $0$ for $\delta$ going to zero. If $i>0$, we can integrate first the variable $r_i$.
Using now that the antiderivative of the function $y/(C+ay^2)$, for a constant $a$, is given by $\log(C+ay^2)/(2a)$ we obtain
\begin{align}\label{eq:integralrec}
&\nonumber C_1\int_{0}^{\ve}  \cdots  \int_{0}^{\ve} \frac{r_{\lceil i \rceil}}{C- \log(\delta)
r_{\lceil i \rceil}^2} \prod_{p=1}^i dr_{p}\\
&\=\frac{C_1}{-2\log(\delta)}\int_{0}^{\ve}  \cdots  \int_{0}^{\ve}
  \log\Bigl(1-\log(\delta)\frac{\ve^2}{C}r_{\lceil i -1\rceil}\Bigr) \prod_{p=1}^{i-1} dr_{p}\,.
\end{align}
If $i=1$, the previous expression goes to zero for $\delta$ going to zero.
If $i\geq 2$, we can integrate again the $r_{i-1}$ variable. Using that the antiderivative
of $\log(1+ay)$ is given by $(1/a+y)\log(1+ay)-y$ we obtain
\begin{align}\label{eq:finalint}
&\nonumber \frac{-C_1}{2\log(\delta)}\int_{0}^{\ve}  \cdots  \int_{0}^{\ve}
\log\Bigl(1-\log(\delta)\frac{\ve^2}{C}r_{\lceil i -1\rceil}\Bigr) \prod_{p=1}^{i-1} dr_{p}\\
&\nonumber\=\frac{-C_1}{2\log(\delta)}\int_{0}^{\ve}  \cdots  \int_{0}^{\ve}
\left(\frac{1}{-\log(\delta)\frac{\ve^2}{C}r_{\lceil i -2\rceil}+\ve}\cdot
\log\Bigl(1-\log(\delta)\frac{\ve^3}{C}r_{\lceil i -2\rceil}\Bigr)-\ve \right) \prod_{p=1}^{i-2} dr_{p}\\
&\,\leq\, \frac{C_1\ve^{i-1}}{2\log(\delta)}+\frac{C_1}{-2\ve\log(\delta)}
\int_{0}^{\ve}  \cdots  \int_{0}^{\ve} \log\Bigl(1-\log(\delta)\frac{\ve^3}{C}r_{\lceil i -2\rceil}\Bigr) \prod_{p=1}^{i-2} dr_{p}
\end{align}
where the last inequality is given by bounding using the bound $r_{\lceil i -2\rceil}\geq 0$ for the denominator of the fraction appearing in front of the logarithm in the second line.
The first term of \eqref{eq:finalint} is clearly going to zero for $\delta$ going to zero.
The second term of \eqref{eq:finalint}  has the same shape as the integral given by \eqref{eq:integralrec}, so by induction we can show that the expression goes to zero for $\delta$ going to zero. 
\par
Hence we have shown that the initial expression given by the integral over $T^{(-i),j}_\delta$ of the product of building blocks involving the $x_{(-i),j}$ variable is going to zero for $\delta$ going to zero, as we wanted.
\end{proof}

%% file: welded_II.tex
%
%
%
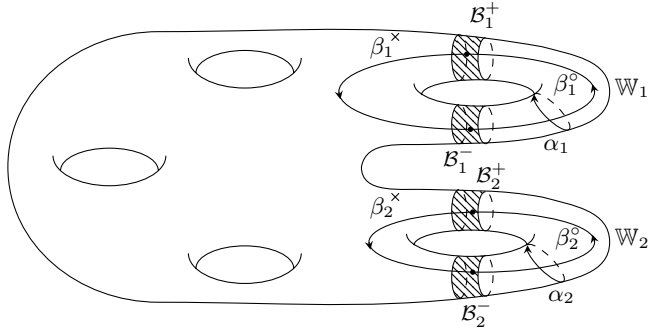
\begin{figure}[htb]
	\centering
\begin{tikzpicture}[>=stealth] 

\begin{scope}
\coordinate (P1) at (0,0);
\coordinate (P2) at (2,1.8);
\coordinate (P3) at (7.1,1.6);
\coordinate (P4) at (8,1);
\coordinate (P5) at (7.1,.4);
\coordinate (P6) at (5.6,.3);
\coordinate (P7) at (4.7,0);
\coordinate (P8) at (5.6,-.3);
\coordinate (P9) at (7.1,-.4);
\coordinate (P10) at (8,-1);
\coordinate (P11) at (7.1,-1.6);
\coordinate (P12) at (2,-1.8);

\path[draw,name path=P1--P1] (P1) to[out=90,in=180] (P2) 
to[out=0,in=170] (P3)
to[out=-15,in=90] (P4)
to[out=270,in=15] (P5)
to[out=190,in=2] (P6)                 
to[out=182,in=80] (P7)
to[out=270,in=178] (P8)
to[out=-2,in=170] (P9)
to[out=-15,in=90] (P10)							      
to[out=270,in=15] (P11)
to[out=190,in=0] (P12)
to[out=180,in=270] (P1);


\draw (6.1,1) ellipse (1.7cm and .5cm);
\draw [->]   (4.4,1) -- (4.42,0.9);
\draw[] (7,.96) arc(0:180:.75cm and .20cm) (7.1,1.11) arc(0:-180:.85cm and .3cm);

\fill[] (6.1,1.5) circle (1pt); 
\draw[] (6,1.15) arc(270:90:.1cm and .305cm);
\draw[dashed] (6,1.15) arc(-90:90:.1cm and .305cm); 
\draw[] (6.35,1.16) arc(270:90:.1cm and .28cm) coordinate (righthightop);
\draw[dashed] (6.35,1.16) arc(-90:90:.1cm and .28cm); 
\fill [pattern=north west lines] 
(6,1.15) arc(270:90:.1cm and .305cm) -- (righthightop)-- (righthightop) arc(-270:-90:.1cm and .28cm); 
\node [above] at (righthightop) {\footnotesize $\cBB_1^+$};

\fill[] (6.15,.5) circle (1pt); 
\draw[] (6,.31) arc(270:90:.1cm and .26cm);
\draw[dashed] (6,.31) arc(-90:90:.1cm and .26cm); 
\draw[] (6.35,.31) arc(270:90:.1cm and .25cm) coordinate (righthighbot);
\draw[dashed] (6.35,.31) arc(-90:90:.1cm and .25cm);
\fill [pattern=north west lines] 
(6,.31) arc(270:90:.1cm and .26cm) -- (righthighbot)-- (righthighbot) arc(-270:-90:.1cm and .25cm); 
\node [below ] at (6,.4) {\footnotesize $\cBB_1^-$};

\begin{scope}[xshift=5.1cm,yshift=1.68cm]
\draw[name path=AB] (0,0)  -- (.1,.1);
\draw[name path=CD] (0,.1) -- (.1,0);
\draw [name intersections={of=AB and CD, by=x}] (x) node [inner sep=0pt, anchor=north east] {\small$\beta_1$};
\end{scope}

\coordinate[](P) at (7.37,.55);
\begin{scope}[rotate around={45:(P)}]
\draw[-{>[flex=0.85]}] (7.4,.45) arc(270:90:.1cm and .33cm);
\draw[dashed] (7.4,.45) arc(-90:90:.1cm and .33cm); 
\end{scope}
\draw (7.3,0.27) node[]{\small $\alpha_1$};
\draw (7.43,1.1) node[]{\small $\beta_1^\circ$}; 
\draw [->]   (7.8,1) -- (7.77,1.1);

\draw[] (2,-.05) arc(0:180:.65cm and .30cm) (2.1,.15) arc(0:-180:.75cm and .4cm);
\begin{scope}[xshift=1.4cm,yshift=-1.3cm]
\draw[] (2.4,-.05) arc(0:180:.65cm and .30cm) (2.5,.15) arc(0:-180:.75cm and .4cm);
\end{scope}
\begin{scope}[xshift=1.4cm,yshift=1.3cm]
\draw[] (2.4,-.05) arc(0:180:.65cm and .30cm) (2.5,.15) arc(0:-180:.75cm and .4cm);
\end{scope}
\end{scope}

\begin{scope}
\draw[] (6.3,-1) ellipse (1.5cm and .4cm);
\draw [->]   (4.8,-1) -- (4.82,-1.1);
\draw[] (6.9,-1.04) arc(0:180:.75cm and .20cm) (7,-.89) arc(0:-180:.85cm and .3cm);

\fill[] (6.18,-.6) circle (1pt); 
\draw[] (6,-.85) arc(270:90:.1cm and .275cm) coordinate (leftlowtop);
\draw[dashed] (6,-.85) arc(-90:90:.1cm and .275cm); 
\draw[] (6.35,-.85) arc(270:90:.1cm and .26cm) coordinate (rightlowtop);
\draw[dashed] (6.35,-.85) arc(-90:90:.1cm and .26cm); 
\fill [pattern=north west lines] 
(6,-.85) arc(270:90:.1cm and .275cm) -- (rightlowtop)-- (rightlowtop) arc(-270:-90:.1cm and .26cm)--(6,-.85); 
\node [below right] at (6.1,.18) {\footnotesize $\cBB_2^+$};

\fill[] (6.18,-1.4) circle (1pt); 
\draw[] (6,-1.75) arc(270:90:.1cm and .28cm) coordinate (leftlowbot);
\draw[dashed] (6,-1.75) arc(-90:90:.1cm and .28cm); 
\draw[] (6.35,-1.72) arc(270:90:.1cm and .27cm) coordinate (rightlowbot);
\draw[dashed] (6.35,-1.72) arc(-90:90:.1cm and .27cm);
\fill [pattern=north west lines] 
(6,-1.75) arc(270:90:.1cm and .28cm) -- (rightlowbot)-- (rightlowbot) arc(-270:-90:.1cm and .27cm)-- (6,-1.75); 
\node [below right] at (5.9,-1.65) {\footnotesize $\cBB_2^-$};
%

\begin{scope}[xshift=5.1cm,yshift=-.45cm]
\draw[name path=AB] (0,0)  -- (.1,.1);
\draw[name path=CD] (0,.1) -- (.1,0);
\draw [name intersections={of=AB and CD, by=x}] (x) node [inner sep=0pt, anchor=north east] {\small $\beta_2$};
\end{scope}
%
\coordinate[](P2) at (7.48,-1.5);
\begin{scope}[rotate around={45:(P2)}]
\draw[-{>[flex=0.85]}] (7.4,-1.45) arc(270:90:.1cm and .34cm);
\draw[dashed] (7.4,-1.45) arc(-90:90:.1cm and .34cm); 
\end{scope}
\draw (7.36,-1.73) node[]{\small$\alpha_2$};
\draw (7.45,-.95) node[]{\small $\beta_2^\circ$};
\draw [->]   (7.8,-1) -- (7.77,-0.9);
\end{scope}

\node [above] at (8.3,-1.2) {\small $\WW_2$};
\node [above] at (8.3,0.8) {\small $\WW_1$};

\end{tikzpicture}
\caption{Decomposing the $\beta_i$ into exterior
and interior of the plumbing fixture}
	\label{cap:betadecomp}
\end{figure}